# Locally complete intersection homomorphisms and a conjecture of Quillen on the vanishing of cotangent homology

By Luchezar L. Avramov*

## Abstract

Classical definitions of locally complete intersection (l.c.i.) homomorphisms of commutative rings are limited to maps that are essentially of finite type, or flat. The concept introduced in this paper is meaningful for homomorphisms $\varphi \colon R \longrightarrow S$ of commutative noetherian rings. It is defined in terms of the structure of $\varphi$ in a formal neighborhood of each point of $\operatorname{Spec} S$. We characterize the l.c.i. property by different conditions on the vanishing of the André-Quillen homology of the $R$-algebra $S$. One of these descriptions establishes a very general form of a conjecture of Quillen that was open even for homomorphisms of finite type: If $S$ has a finite resolution by flat $R$-modules and the cotangent complex $\operatorname{L}(S|R)$ is quasi-isomorphic to a bounded complex of flat $S$-modules, then $\varphi$ is l.c.i. The proof uses a mixture of methods from commutative algebra, differential graded homological algebra, and homotopy theory. The l.c.i. property is shown to be stable under a variety of operations, including composition, decomposition, flat base change, localization, and completion. The present framework allows for the results to be stated in proper generality; many of them are new even with classical assumptions. For instance, the stability of l.c.i. homomorphisms under decomposition settles an open case in Fulton's treatment of orientations of morphisms of schemes.

## Table of Contents

*Research partially supported by a grant from the National Science Foundation.
1991 *Mathematics Subject Classification*. 13D03, 13H10, 14E40, 14M10.

**Introduction**

The concept of regularity—of distinctly geometric origin—has taken a central place in the study of commutative noetherian rings. This is in part due to the existence of a Cohen presentation of each complete local ring as a homomorphic image of a regular ring.

From a homological perspective, a local ring is closest to being regular if it is *complete intersection* (or: c.i.) in the sense that the defining ideal of some Cohen presentation of its completion is generated by a regular sequence. A noetherian ring is *locally complete intersection* (or: l.c.i.) if its localizations at all prime ideals are complete intersections.

The relative version of the notion of regularity is well established: a homomorphism $\varphi \colon R \longrightarrow S$ is *regular* if it is flat and has geometrically regular fibers. Foundational work of Grothendieck [26], Lichtenbaum and Schlessinger [32], André [1], [3], and Quillen [37], characterized regularity by the vanishing for $n \geq 1$ of the functors $\mathrm{D}_n(S\,|\,R, -)$ of André-Quillen (or: cotangent) homology.

In contrast, no general notion of l.c.i. homomorphism has emerged. For philosophical, historical, and practical reasons, such a concept has to accommodate the following cases:

- When $R$ is regular, $\varphi$ is l.c.i. precisely when $S$ is l.c.i.
- When $\varphi$ is flat, it is l.c.i. if and only if all its nontrivial fiber rings are l.c.i.
- When $\varphi$ can be factored as a regular map followed by a surjection $\varphi'$, it is l.c.i. if and only if in some factorization $\operatorname{Ker} \varphi'$ can be locally generated by a regular sequence.

To reconcile these notions when several apply, it was proved in *loc. cit.* that each one is equivalent to the vanishing of $\mathrm{D}_n(S\,|\,R, -)$ for $n \geq 2$. Maps with this property were called weakly c.i. by Illusie [31] who remarked that useful properties, like transitivity or flat base change, follow directly from formal properties of cotangent homology. While settling functorial questions, the homological definition gave no approach to structural properties; for instance, it was not known if a weakly c.i. map is locally of finite flat dimension, that is, if $S_{\mathfrak{q}}$ has a finite resolution by flat $R$-modules for each prime ideal $\mathfrak{q} \subseteq S$.

For homomorphisms of *noetherian rings* we introduce an l.c.i. notion locally, by using a relative version of Cohen structure theory developed jointly with Foxby and B. Herzog [14]. We show that these maps coincide with the weakly c.i. homomorphisms, then characterize them by the vanishing of $\mathrm{D}_2(S\,|\,R, -)$, and by various other vanishing conditions.

In particular, we prove that $\varphi$ is l.c.i. if and only if $\varphi$ is locally of finite flat dimension and $\mathrm{D}_n(S\,|\,R, -) = 0$ for all $n \gg 0$; for maps essentially of finite type the 'only if' part is well known, and the converse was conjectured by Quillen [37]. When $\mathbb{Q} \subseteq S$ and $\varphi$ is locally of finite flat dimension, we proved in joint work with Halperin [16] that if $\varphi$ is not l.c.i. then $\mathrm{D}_n(S \,|\, R, -) \neq 0$ for $n \gg 0$; this is strengthened below: the vanishing of $\mathrm{D}_n(S \,|\, R, -)$ for *any* single $n \geq 1$ implies that $\varphi$ is l.c.i.

An interplay of structural and homological arguments allows for a comprehensive study of l.c.i. maps, in the framework of a joint program with Foxby [11], [12], [13] to classify ring homomorphisms according to their local structure. We establish the stability of the new class under composition, decomposition, flat base change, localization, completion, and clarify its role in the transfer of l.c.i. properties between source and target rings.

An application deals with the functoriality of canonical orientations assigned to morphisms of schemes that are l.c.i., or flat. Fulton [24] proves that orientations of such morphisms $f \colon X \longrightarrow Y$ and $g \colon Y \longrightarrow Z$ satisfy $[gf] = [f][g]$ in five out of the six possible cases. In the open case when $f$ is flat and both $g$ and $gf$ are l.c.i. our decomposition theorem for l.c.i. homomorphisms shows that $f$ is l.c.i., so the formula holds as well.

Homomorphisms of commutative rings and their simplicially defined homology theory are the subject of this investigation, so a significant role is predictably played by commutative algebra and homotopy theory. The gap is bridged by DG (= differential graded) homological algebra, viewed alternatively as an extension of the former and a linearization of the latter. It produces invariants of local homomorphisms that are often more computable than those defined in terms of cotangent complexes.

The results that follow have evolved over a long period of time, through several preliminary versions and oral expositions[1]. As a consequence, the paper has benefited from direct or indirect input of several people. My thinking on l.c.i. homomorphisms has been influenced by two collaborations: with Steve Halperin on links between local algebra and rational homotopy, and with Hans-Bjørn Foxby on local properties of ring homomorphisms. Javier Majadas showed me a substantial shortcut in the proof of Lemma (1.7). Srikanth Iyengar pointed out that an argument I had developed for a different purpose could be used in the proof of Theorem (3.4). Haynes Miller made me aware of Umeda's work [34] and kept reminding me to produce a final version. I am grateful to all of them.

---

[1] A discussion of properties of l.c.i. homomorphisms and the structure of the proof of Quillen's Conjecture is contained in the extended abstract of my talk "Locally complete intersection homomorphisms and vanishing of André-Quillen homology" [Commutative algebra. International conference, Vechta, 1994 (W. Bruns, J. Herzog, M. Hochster, U. Vetter, eds.), Runge, Cloppenburg, 1994, pp. 20–24].

## 1. Vanishing theorems

Throughout this section $\varphi\colon R \to S$ is a homomorphism of noetherian rings.

The André-Quillen homology $\mathrm{D}_n(S\,|\,R,\,N)$ of the $R$-algebra $S$ with coefficients in an $S$-module $N$ is the $n^{\mathrm{th}}$ homology module of $\mathrm{L}(S\,|\,R) \otimes_S N$, where $\mathrm{L}(S\,|\,R)$ is the *cotangent complex* of $\varphi$, uniquely defined in the derived category of $S$-modules $\mathbf{D}(S)$ (cf. [3], [37]). The vanishing of $\mathrm{D}_n(S\,|\,R,\,-)$ for all $n > m$ means that $\mathrm{L}(S\,|\,R)$ is isomorphic in $\mathbf{D}(S)$ to a bounded complex $F$ of flat $S$-modules with $F_n = 0$ for $n > m$; we then say that it has *flat dimension* at most $m$, and write $\mathrm{fd}_S\,\mathrm{L}(S\,|\,R) \leq m$.

Regularity is characterized by the vanishing of the first cotangent homology functor. We quote that benchmark result in the definitive version of André [3, (S.30)].

(1.1) FIRST VANISHING THEOREM. *The following conditions are equivalent*:

 (i) $\varphi$ *is regular*.
 (ii) $\mathrm{D}_1(S\,|\,R,\,-) = 0$.
 (iii) $\mathrm{fd}_S\,\mathrm{L}(S\,|\,R) = 0$.

For the special types of homomorphisms reviewed in the introduction, the vanishing of $\mathrm{D}_2(S\,|\,R,\,-)$ is classically known to be equivalent to the corresponding l.c.i. notion. To describe the vanishing condition in general, we recall a structure theorem.

If the homomorphism $\varphi$ is *local*, in the sense that both rings are local and $\varphi(\mathfrak{m}) \subseteq \mathfrak{n}$ where $\mathfrak{m}$ is the unique maximal ideal of $R$ and $\mathfrak{n}$ is that of $S$, then by [14, (1.1)] the composition $\dot{\varphi}\colon R \to \widehat{S}$ of $\varphi$ with the completion map $S \to \widehat{S}$ appears in a commutative diagram of local homomorphisms of local rings

$$\begin{array}{ccc} & R' & \\ {\scriptstyle \dot\varphi}\nearrow & & \searrow{\scriptstyle \varphi'} \\ R & \xrightarrow[\dot\varphi]{} & \widehat{S} \end{array}$$

where $\dot\varphi$ is flat, $\varphi'$ is surjective, the ring $R'$ is complete, and the ring $R'/\mathfrak{m}R'$ is regular. Such a diagram is called a *Cohen factorization* of $\dot\varphi$.

*Definition.* A local homomorphism $\varphi\colon R \to S$ is *complete intersection*, or *c.i.*, at $\mathfrak{n}$, if in some Cohen factorization $\mathrm{Ker}\,\varphi'$ is generated by a regular sequence (this property does not depend on the choice of Cohen factorization; cf. (3.3)).

A homomorphism of noetherian rings $\varphi \colon R \to S$ is c.i. at a prime ideal $\mathfrak{q}$ of $S$ if the induced local homomorphism $\varphi_{\mathfrak{q}} \colon R_{\mathfrak{q} \cap R} \to S_{\mathfrak{q}}$ is c.i. at $\mathfrak{q} S_{\mathfrak{q}}$. A homomorphism that has this property at all $\mathfrak{q} \in \operatorname{Spec} S$ is said to be *locally complete intersection*, or *l.c.i.*

With this notion, we have:

(1.2) SECOND VANISHING THEOREM. *The following conditions are equivalent*:

 (i) $\varphi$ *is locally complete intersection.*
 (ii) $\operatorname{D}_2(S \,|\, R, -) = 0$.
(iii) $\operatorname{fd}_S \operatorname{L}(S \,|\, R) \leq 1$.

The proof, given at the end of this section, uses the existence of Cohen factorizations and only standard properties of André-Quillen homology. For expository reasons, we continue with a discussion of vanishing results proved at the end of Section 4.

Let $\operatorname{fd}_R M$ denote the flat dimension (also called Tor-dimension) of an $R$-module $M$. We say that $\varphi$ is *locally of finite flat dimension* if $\operatorname{fd}_R S_{\mathfrak{q}}$ is finite for all $\mathfrak{q} \in \operatorname{Spec} S$. This condition is clearly implied by the finiteness of $\operatorname{fd}_R S$, and is equivalent to it in many cases, e.g. when $R$ has finite Krull dimension; cf. [5].

For maps essentially of finite type the *only if* part of the next theorem is well known. Quillen [37, (5.7)] conjectured that the converse holds as well. This was proved by Avramov and Halperin [16] in characteristic zero. We establish a very general form of

(1.3) QUILLEN'S CONJECTURE. *The homomorphism $\varphi$ is l.c.i. if and only if it is locally of finite flat dimension and $\operatorname{fd}_S \operatorname{L}(S \,|\, R)$ is finite.*

The next result represents a partial strengthening of Quillen's conjecture. Note that the condition on $m$ poses no restriction when $n = 3$, or when $S$ is an algebra over $\mathbb{Q}$.

(1.4) RIGIDITY THEOREM. *Let $m \geq 2$ be an integer, such that $(m-1)!$ is invertible in $S$.*

*If $\varphi$ is locally of finite flat dimension and $\operatorname{D}_n(S \,|\, R, -) = 0$ for some $n$ with $3 \leq n \leq 2m - 1$, then $\varphi$ is locally complete intersection.*

In view of the preceding results we extend to all homomorphisms of noetherian rings another conjecture, proposed by Quillen [37, (5.6)] for maps essentially of finite type:

CONJECTURE. *If $\operatorname{fd}_S \operatorname{L}(S \,|\, R) < \infty$, then $\operatorname{fd}_S \operatorname{L}(S \,|\, R) \leq 2$.*

We are able to verify it when one of the rings $R$ or $S$ is l.c.i.; for a more precise statement, we use the map ${}^a\varphi \colon \operatorname{Spec} S \to \operatorname{Spec} R$ induced by $\varphi$.

(1.5) THEOREM ON L.C.I. RINGS. *If $R$ is c.i. on the image of ${}^a\varphi$ and $\operatorname{fd}_S \operatorname{L}(S|R) < \infty$, then $S$ is l.c.i.*

*If $S$ is l.c.i. and $\operatorname{fd}_S \operatorname{L}(S|R) < \infty$, then $R$ is c.i. on the image of ${}^a\varphi$.*

*If $R$ is c.i. on the image of ${}^a\varphi$ and $S$ is l.c.i., then $\operatorname{fd}_S \operatorname{L}(S|R) \leq 2$.*

To prepare for the proof of (1.2) we recall a few basic results on cotangent homology.

(1.6) *Remarks.* Let $(R', \mathfrak{m}', \ell)$ be a local ring.

(1) For an ideal $\mathfrak{a} \subsetneq R'$ the following are equivalent: (i) $\mathfrak{a}$ is generated by a regular sequence; (ii) $\operatorname{D}_2(R'/\mathfrak{a} \,|\, R', \ell) = 0$; (iii) $\operatorname{D}_n(R'/\mathfrak{a} \,|\, R', \ell) = 0$ for $n \geq 2$ (cf. [3, (6.25)]).

(2) The ring $R'$ is regular if and only if $\operatorname{D}_2(\ell \,|\, R', \ell) = 0$ (cf. [3, (6.26)]).

(3) If $k \to \ell$ is a field extension, then $\operatorname{fd}_\ell \operatorname{L}(\ell|k) \leq 1$ (cf. [3, (7.4)]).

(1.7) LEMMA. *If $\varphi \colon (R, \mathfrak{m}, k) \to (S, \mathfrak{n}, \ell)$ is a local homomorphism, $\sigma \colon S \to \widehat{S}$ is the completion map, and $R \xrightarrow{\dot{\varphi}} R' \xrightarrow{\varphi'} \widehat{S}$ is a Cohen factorization of $\dot{\varphi}$, then the canonical map $\operatorname{D}_n(\sigma \,|\, \dot{\varphi}, \ell) \colon \operatorname{D}_n(S \,|\, R, \ell) \to \operatorname{D}_n(\widehat{S} \,|\, R', \ell)$ is an isomorphism for $n \geq 2$.*

*Proof.* For each $n$ the map in question is equal to $\operatorname{D}_n(\widehat{S}\,|\,\dot{\varphi}, \ell) \circ \operatorname{D}_n(\sigma\,|\,R, \ell)$.

Flat base change gives $\operatorname{D}_n(\widehat{S}\,|\,S, \ell) \cong \operatorname{D}_n(\ell\,|\,\ell, \ell) = 0$ for all $n$ so the Jacobi-Zariski exact sequence of the decomposition $\dot{\varphi} = \sigma\varphi$ shows that $\operatorname{D}_n(\sigma\,|\,R, \ell)$ is an isomorphism for each $n$. Thus, it suffices to prove that $\operatorname{D}_n(\widehat{S}\,|\,\dot{\varphi}, \ell)$ is bijective for $n \geq 2$.

Set $\overline{R}' = R'/\mathfrak{m}R'$. From the Jacobi-Zariski exact sequence of $k \to \overline{R}' \to \ell$,

$$\operatorname{D}_{n+1}(\ell\,|\,\overline{R}', \ell) \to \operatorname{D}_n(\overline{R}'\,|\,k, \ell) \to \operatorname{D}_n(\ell\,|\,k, \ell).$$

Since $\operatorname{D}_n(\ell\,|\,k, \ell) = 0$ for $n \geq 2$ by (1.6.3) and $\operatorname{D}_{n+1}(\ell\,|\,\overline{R}', \ell) = 0$ for $n \geq 1$ by (1.6.2), we get $\operatorname{D}_n(\overline{R}'\,|\,k, \ell) = 0$ for $n \geq 2$. Flat base change (cf. [3, (4.54)]) yields isomorphisms $\gamma_n \colon \operatorname{D}_n(R'\,|\,R, \ell) \cong \operatorname{D}_n(\overline{R}'\,|\,k, \ell)$ for all $n \in \mathbb{Z}$, so using the Jacobi-Zariski exact sequence

$$\operatorname{D}_n(R'\,|\,R,\,\ell) \xrightarrow{\operatorname{D}_n(\varphi'\,|\,R,\ell)} \operatorname{D}_n(\widehat{S}\,|\,R,\,\ell) \xrightarrow{\operatorname{D}_n(\widehat{S}\,|\,\dot{\varphi},\ell)} \operatorname{D}_n(\widehat{S}\,|\,R',\,\ell)$$
$$\xrightarrow{\eth_n} \operatorname{D}_{n-1}(R'\,|\,R,\,\ell)$$

of the Cohen factorization $\dot{\varphi} = \varphi'\dot{\varphi}$ we conclude that $\operatorname{D}_n(\widehat{S}\,|\,\dot{\varphi}, \ell)$ is bijective for $n \geq 3$ and injective for $n = 2$. Comparison of another segment of that sequence

with the one for the homomorphisms $k \to \overline{R}' \to \ell$ yields a commutative diagram with exact rows

$$\begin{array}{ccccccc}
D_2(\widehat{S}\,|\,R,\,\ell) & \xrightarrow{D_2(\widehat{S}\,|\,\dot{\varphi},\,\ell)} & D_2(\widehat{S}\,|\,R',\,\ell) & \xrightarrow{\eth_2} & D_1(R'\,|\,R,\,\ell) & \xrightarrow{D_1(\varphi'\,|\,R,\,\ell)} & D_1(\widehat{S}\,|\,R,\,\ell) \\
\downarrow & & & & \gamma_1 \downarrow \cong & & \downarrow \\
D_2(\ell\,|\,\overline{R}',\,\ell) & & \longrightarrow & & D_1(\overline{R}'\,|\,k,\,\ell) & \longrightarrow & D_1(\ell\,|\,k,\,\ell).
\end{array}$$

It implies that $\eth_2 = 0$, so the map $D_2(\widehat{S}\,|\,\dot{\varphi},\,\ell)$ is surjective. □

Combining the lemma with (1.6.1), we get:

(1.8) PROPOSITION 1.1. *A local homomorphism $R \to (S, \mathfrak{n}, \ell)$ is complete intersection at $\mathfrak{n}$ if and only if $D_2(S\,|\,R,\,\ell) = 0$; when this is the case, $D_n(S\,|\,R,\,\ell) = 0$ for $n \geq 2$.*

Recall that André-Quillen homology localizes perfectly.

(1.9) *Remark.* Let $\varphi \colon R \to S$ be a homomorphism of commutative rings. For each $n \in Z$ and each $\mathfrak{q} \in \mathrm{Spec}\, S$ there is a isomorphism $D_n(R|S, -)_\mathfrak{q} \cong D_n(S_\mathfrak{q}|R_{\mathfrak{q} \cap S}, -_\mathfrak{q})$ of functors on the category of $S$-modules, (cf. [3, (4.59) and (5.27)]).

*Proof of Theorem* (1.2). By (1.8) and (1.9), the following conditions are equivalent:

(i)  $\varphi$ is l.c.i.;
(ii') $D_2(S\,|\,R,\,k(\mathfrak{q})) = 0$ for each $\mathfrak{q} \in \mathrm{Spec}\, S$;
(iii') $D_n(S\,|\,R,\,k(\mathfrak{q})) = 0$ for each $\mathfrak{q} \in \mathrm{Spec}\, S$ and all $n \geq 2$.

If $D_2(S\,|\,R,\,-) = 0$, then (ii') holds by (1.9).
On the other hand, if (iii') holds, then $D_n(S\,|\,R,\,-) = 0$ for $n \geq 2$ by [3, (S.29)]. □

## 2. Eilenberg-Zilber quasi-isomorphisms

The classical Eilenberg-Zilber theorem shows that the normalized chain complex of a tensor product of simplicial abelian groups is homotopy equivalent to the tensor product of their normalized chain complexes. We produce a substitute for simplicial modules over a simplicial ring. On the way, we introduce some notation for DG and simplicial algebra.

(2.1) *DG algebra* (cf. [34], [27], [15]). DG objects have differentials of degree $-1$, denoted ubiquitously $\partial$. Morphisms of DG objects are chain maps of degree 0 that preserve the appropriate structure; quasi-isomorphisms are morphisms that induce isomorphism in homology. The functor $(-)^\natural$ forgets differentials; $|x|$ denotes the degree of an element $x$. Graded algebras are trivial in negative degrees; graded modules are bounded below.

Elementary proofs of the next two assertions can be found in [10, §1.3].

*Remark.* Let $M$ be a right DG module over $A$, such that $M^\natural$ is a free $A^\natural$-module.

(1) If $\nu\colon N' \to N$ is a quasi-isomorphism of left DG modules, then the induced map $M \otimes_A \nu\colon M \otimes_A N' \to M \otimes_A N$ is a quasi-isomorphism.

(2) If $\mu\colon M' \to M$ is a quasi-isomorphism of right DG modules over $A$, and $M'^\natural$ is free over $A^\natural$, then for each left DG module $N$ the induced map $\mu \otimes_A N\colon M' \otimes_A N \to M \otimes_A N$ is a quasi-isomorphism.

(2.2) *Simplicial algebra* (cf. [20], [37]). The face operators $d_i^n\colon \mathcal{G}_n \to \mathcal{G}_{n-1}$ ($i = 0, \ldots, n$) of a simplicial object $\mathcal{G}$ are morphisms in the corresponding category. The normalization functor N from simplicial abelian groups to nonnegatively graded chain complexes has $(\mathrm{N}\mathcal{G})_n = \bigcap_{i=1}^n \mathrm{Ker}\,(d_i^n)$ and $\partial_n\colon (\mathrm{N}\mathcal{G})_n \to (\mathrm{N}\mathcal{G})_{n-1}$ equal to the restriction of $d_0^n$. The homotopy of $\mathcal{G}$ is the graded abelian group $\pi(\mathcal{G}) = \mathrm{H}(\mathrm{N}\mathcal{G})$. A weak equivalence is a homomorphism of simplicial groups that induces an isomorphism in homotopy. Normalization has a quasi-inverse, given by the Dold-Kan functor K (cf. [20, (3.6)]). The functors N and K transform weak equivalences and quasi-isomorphisms into each other.

The functor that assigns to a bisimplicial object its diagonal simplicial subobject is denoted $\Delta$. Thus, if $\mathcal{B}$ is a bisimplicial abelian group, then $\Delta\mathcal{B}$ has $(\Delta\mathcal{B})_n = \mathcal{B}_{n,n}$ and face operators $d_{i,i}^n$ for $i = 0.\ldots,n$. We denote $\mathrm{N}^\mathrm{v}$ (respectively, $\mathrm{N}^\mathrm{h}$) the normalization functor applied to the columns (respectively, rows) of $\mathcal{B}$. By the Eilenberg-Zilber-Cartier theorem [20, (2.9)] there is a natural weak equivalence $\mathrm{Tot}\,(\mathrm{N}^\mathrm{v}\mathrm{N}^\mathrm{h}\mathcal{B}) \to \mathrm{N}\Delta\mathcal{B}$, where $\mathrm{Tot}\,(B_\bullet)$ denotes the total complex associated to a double complex $B_\bullet$.

Let $\mathcal{M}$ be a simplicial right module; $\mathcal{M}$ is *cofibrant* if for each surjective weak equivalence $\alpha\colon \mathcal{L}' \to \mathcal{L}$ of right simplicial $\mathcal{A}$-modules and each homomorphism $\gamma\colon \mathcal{M} \to \mathcal{L}$ there is a homomorphism $\beta\colon \mathcal{M} \to \mathcal{L}'$ such that $\alpha\beta = \gamma$. If $\mathcal{N}$ is a simplicial left $\mathcal{A}$-module, then $\mathcal{M}\bar{\otimes}_\mathcal{A}\mathcal{N}$ is the simplicial abelian group with $(\mathcal{M}\bar{\otimes}_\mathcal{A}\mathcal{N})_n = \mathcal{M}_n \otimes_{\mathcal{A}_n} \mathcal{N}_n$ and diagonal simplicial operators. Shuffle products give the normalization $N\mathcal{A}$ a structure of DG algebra, and make $N\mathcal{M}$ and $N\mathcal{N}$ into a right and left DG module over $N\mathcal{A}$, respectively.

PROPOSITION. *Let $\mathcal{N}$ be a simplicial left module over a simplicial ring $\mathcal{A}$, let $\mathcal{M}$ be a cofibrant simplicial right module over $\mathcal{A}$, and let $M'$ be a right DG module over the DG ring $N\mathcal{A}$, such that $M'^{\natural}$ is free over $N\mathcal{A}^{\natural}$.*

*If $\mu \colon M' \to N\mathcal{M}$ is a quasi-isomorphism of right DG modules over $N\mathcal{A}$, then the composition of $\mu \otimes_{N\mathcal{A}} N\mathcal{N} \colon M' \otimes_{N\mathcal{A}} N\mathcal{N} \to N\mathcal{M} \otimes_{N\mathcal{A}} N\mathcal{N}$ with the canonical map $N\mathcal{M} \otimes_{N\mathcal{A}} N\mathcal{N} \to N(\mathcal{M} \bar{\otimes}_{\mathcal{A}} \mathcal{N})$ is a quasi-isomorphism.*

*Proof.* Illusie [31, (3.3.3.8)] constructs an exact sequence of simplicial right $\mathcal{A}$-modules

$$\mathcal{P}_{\bullet}^+ \colon \ldots \to \mathcal{P}_{[p]} \xrightarrow{\delta_{[p]}} \mathcal{P}_{[p-1]} \to \ldots \to \mathcal{P}_{[0]} \xrightarrow{\delta_{[0]}} \mathcal{M} \to 0$$

where for each $p \geq 0$ the simplicial $\mathcal{A}$-module $\mathcal{P}_{[p]}$ is equal to $\mathcal{L}_{[p]} \bar{\otimes}_{\mathbb{Z}} \mathcal{A}$, with $\mathcal{L}_{[p]}$ a simplicial abelian group such that $\mathcal{L}_{[p]}$ and $\pi(\mathcal{L}_{[p]})$ are free graded abelian groups.

Let $K^h$ denote the Dold-Kan functor applied to a nonnegative complex of simplicial right $\mathcal{A}$-modules. It produces a bisimplicial group with $p^{\text{th}}$ column a simplicial right module over the simplicial ring $\mathcal{A}$ for each $p$, and $q^{\text{th}}$ row a simplicial right module over the ring $\mathcal{A}_q$ for each $q$; the diagonal is naturally a simplicial right $\mathcal{A}$-module.

Let $\mathcal{M}_{\bullet}$ be the complex of simplicial $\mathcal{A}$-modules with $\mathcal{M}_{[0]} = \mathcal{M}$ and $\mathcal{M}_{[p]} = 0$ for $p \neq 0$. The map $\delta_{[0]}$ defines a morphism $\mathcal{P}_{\bullet} = \mathcal{P}_{\bullet}^+/\mathcal{M} \to \mathcal{M}_{\bullet}$ of complexes of simplicial $\mathcal{A}$-modules. It induces an isomorphism in $\delta$-homology, hence a weak equivalence of bisimplicial $\mathcal{A}$-modules $K^h \mathcal{P}_{\bullet} \to K^h \mathcal{M}_{\bullet}$, and finally a weak equivalence of simplicial right $\mathcal{A}$-modules $\epsilon \colon \mathcal{Q} = \Delta K^h \mathcal{P}_{\bullet} \to \Delta K^h \mathcal{M}_{\bullet} = \mathcal{M}$. This yields a commutative square of homomorphisms of simplicial abelian groups

$$\begin{array}{ccc} \mathcal{Q} \bar{\otimes}_{\mathcal{A}}^{\mathbf{L}} \mathcal{N} & \xrightarrow{\epsilon \bar{\otimes}_{\mathcal{A}}^{\mathbf{L}} \mathcal{N}} & \mathcal{M} \bar{\otimes}_{\mathcal{A}}^{\mathbf{L}} \mathcal{N} \\ \downarrow & & \downarrow \\ \mathcal{Q} \bar{\otimes}_{\mathcal{A}} \mathcal{N} & \xrightarrow{\epsilon \bar{\otimes}_{\mathcal{A}} \mathcal{N}} & \mathcal{M} \bar{\otimes}_{\mathcal{A}} \mathcal{N} \end{array}$$

where $-\bar{\otimes}_{\mathcal{A}}^{\mathbf{L}} -$ denotes the derived tensor product of Quillen [37, II.6] and the vertical arrows are canonical homomorphisms. The top arrow is a weak equivalence along with $\epsilon$. By construction, $\mathcal{Q}_n$ is a free module over $\mathcal{A}_n$ for each $n \geq 0$; by hypothesis, $\mathcal{M}$ is cofibrant, so $\mathcal{M}_n$ is a direct summand of a free $\mathcal{A}_n$-module for each $n \geq 0$. By [37, p. II.6.10] both vertical maps are weak equivalences, hence $\epsilon \bar{\otimes}_{\mathcal{A}} \mathcal{N}$ is a weak equivalence.

Set $A = \mathrm{N}\mathcal{A}$, $M = \mathrm{N}\mathcal{M}$, and $N = \mathrm{N}\mathcal{N}$. The complex $P_\bullet = \mathrm{N}^{\mathrm{v}}\mathcal{P}_\bullet$ of right DG modules over $A$, and the right DG modules $P = \mathrm{Tot}\,(P_\bullet)$ and $Q = \mathrm{N}\mathcal{Q}$, appear in a diagram

$$\begin{array}{ccccc}
P \otimes_A N & = & \mathrm{Tot}\,(\mathrm{N}^{\mathrm{v}}(\mathcal{P}_\bullet) \otimes_A N) & \xrightarrow{\alpha} & \mathrm{Tot}\,(\mathrm{N}^{\mathrm{v}}(\mathcal{P}_\bullet \bar{\otimes}_\mathcal{A} \mathcal{N})) \\
\| & & \downarrow{\scriptstyle \mathrm{Tot}\,((\mathrm{N}^{\mathrm{v}}\beta)\otimes_A N)} & & \downarrow{\scriptstyle \mathrm{Tot}\,(\mathrm{N}^{\mathrm{v}}\beta')} \\
\mathrm{Tot}\,(\mathrm{N}^{\mathrm{v}}\mathrm{N}^{\mathrm{h}}\mathrm{K}^{\mathrm{h}}\mathcal{P}_\bullet) \otimes_A N & = & \mathrm{Tot}\,(\mathrm{N}^{\mathrm{v}}(\mathrm{N}^{\mathrm{h}}\mathrm{K}^{\mathrm{h}}\mathcal{P}_\bullet) \otimes_A N) & \xrightarrow{\alpha'} & \mathrm{Tot}\,(\mathrm{N}^{\mathrm{v}}\mathrm{N}^{\mathrm{h}}\mathrm{K}^{\mathrm{h}}(\mathcal{P}_\bullet \bar{\otimes}_\mathcal{A} \mathcal{N})) \\
\downarrow{\scriptstyle \gamma \otimes_A N} & & & & \downarrow{\scriptstyle \gamma'} \\
\mathrm{N}(\Delta\mathrm{K}^{\mathrm{h}}\mathcal{P}_\bullet) \otimes_A N & \xrightarrow{\eta'} & \mathrm{N}((\Delta\mathrm{K}^{\mathrm{h}}\mathcal{P}_\bullet)\bar{\otimes}_\mathcal{A}\mathcal{N}) & = & \mathrm{N}(\Delta\mathrm{K}^{\mathrm{h}}\mathcal{P}_\bullet \bar{\otimes}_\mathcal{A} \mathcal{N}) \\
\| & & \| & & \\
Q \otimes_A N & \xrightarrow{\eta} & \mathrm{N}(\mathcal{Q}\bar{\otimes}_\mathcal{A}\mathcal{N}) & &
\end{array}$$

of morphisms of chain complexes defined as follows:

- the equalities are canonical identifications;
- $\alpha$ and $\alpha'$ are totalings of shuffle products;
- $\beta\colon \mathcal{P}_\bullet \to \mathrm{N}^{\mathrm{h}}\mathrm{K}^{\mathrm{h}}\mathcal{P}_\bullet$ and $\beta'\colon \mathcal{P}_\bullet\bar{\otimes}_\mathcal{A}\mathcal{N} \to \mathrm{N}^{\mathrm{h}}\mathrm{K}^{\mathrm{h}}\mathcal{P}_\bullet\bar{\otimes}_\mathcal{A}\mathcal{N}$ are Dold-Kan isomorphisms;
- $\gamma$ and $\gamma'$ are Eilenberg-Zilber-Cartier homotopy equivalences;
- $\eta$ and $\eta'$ are shuffle products.

Thus, all vertical maps are quasi-isomorphisms, and the diagram commutes due to the naturality of all the maps involved.

Filtering the chain complexes $\mathrm{Tot}\,(\mathrm{N}^{\mathrm{v}}(\mathcal{P}_\bullet) \otimes_A N)$ and $\mathrm{Tot}\,(\mathrm{N}^{\mathrm{v}}(\mathcal{P}_\bullet\bar{\otimes}_\mathcal{A}\mathcal{N}))$ by the resolution degree of $\mathcal{P}_\bullet$ we get a homomorphism of spectral sequences ${}^r\alpha_{*,*}\colon {}^r\mathrm{E}'_{*,*} \to {}^r\mathrm{E}^{\mathrm{Q}}_{*,*}$ for $r \geq 0$. The map ${}^0\alpha_{p,*}$ appears for each $p$ in a commutative diagram of chain maps

$$\begin{array}{ccccc}
\mathrm{Tot}\,(L_{[p]} \otimes_\mathbb{Z} A) \otimes_A N & \xrightarrow{\zeta \otimes_\mathbb{Z} N} & \mathrm{N}(\mathcal{L}_{[p]}\bar{\otimes}_\mathbb{Z}\mathcal{A}) \otimes_A N & = & P_{[p]} \otimes_A N \\
\| & & & & \downarrow{\scriptstyle {}^0\alpha_{p,*}} \\
\mathrm{Tot}\,(L_{[p]} \otimes_\mathbb{Z} N) & \xrightarrow{\zeta'} & \mathrm{N}(\mathcal{L}_{[p]}\bar{\otimes}_\mathbb{Z}\mathcal{N}) & = & \mathrm{N}(\mathcal{P}_{[p]}\bar{\otimes}_\mathcal{A}\mathcal{N})
\end{array}$$

where $P_{[p]} = \mathrm{N}\mathcal{P}_{[p]}$ and $L_{[p]} = \mathrm{N}\mathcal{L}_{[p]}$, while $\zeta$ and $\zeta'$ are classical Eilenberg-Zilber homotopy equivalences. Thus, ${}^1\alpha_{p,*}$ is bijective and so $\mathrm{H}(\alpha)$ is an isomorphism. We have shown that all the maps in the first diagram are quasi-isomorphisms.

An Eilenberg-Moore resolution of $D$ is a complex of morphisms of right DG modules

$$D_\bullet^+ : \ldots \to D_{[p]} \xrightarrow{\delta_{[p]}} D_{[p-1]} \to \ldots \to D_{[0]} \xrightarrow{\delta_{[0]}} M' \to 0$$

such that the functors $(-)^\natural$ and $\mathrm{H}(-)$, respectively forgetting the internal differentials $\partial$ and computing their homology, yield exact sequences of free graded modules over $A^\natural$ and $\mathrm{H}(A)$, respectively; we refer to [34], [27], [15] for the construction of such resolutions and of a morphism of complexes $\xi_\bullet^+ \colon D_\bullet^+ \to P_\bullet^+$ with $\xi_{-1}^+ = \mu \colon M' \to M$.

Set $D_\bullet = D_\bullet^+/M'$ and $D = \mathrm{Tot}\, D_\bullet$. The morphism $\xi_\bullet \colon D_\bullet = D_\bullet^+/M' \to P_\bullet$ of complexes of right DG modules over $A$ induces a morphism of DG modules $\xi = \mathrm{Tot}\, \xi_\bullet \colon D \to P$, and so a chain map $\xi \otimes_A N \colon D \otimes_A N \to P \otimes_A N$ that respects the filtrations by resolution degree. As a result we get a homomorphism of spectral sequences ${}^r\xi \colon {}^r\mathrm{E}^{\mathrm{EM}} \to {}^r\mathrm{E}'$ for $r \geq 0$ that converges to $\mathrm{H}(\xi \otimes_A N) \colon \mathrm{H}(D \otimes_A N) \to \mathrm{H}(P \otimes_A N)$. By construction, the map $\mathrm{H}(\xi) \colon \mathrm{H}(D_\bullet) \to \mathrm{H}(P_\bullet)$ is a morphism of free resolutions over $\mathrm{H}(A)$ and induces an isomorphism $\mathrm{H}(\mu)$. We conclude that ${}^2\xi$ is an isomorphism, hence so is $\mathrm{H}(\xi \otimes_A N)$.

Let $\epsilon' \colon D \to M'$ be the quasi-isomorphism of right DG modules induced by $\delta_{[0]}$. As the $A^\natural$-modules $D^\natural$ and $M'^\natural$ are free, respectively by construction and by hypothesis, $\epsilon' \otimes_A N$ is a quasi-isomorphism by (2.1.2). Thus, we now have a commutative diagram of chain complexes in which all arrows adorned by $\simeq$ are quasi-isomorphisms

$$\begin{array}{ccccccc}
D \otimes_A N & \xrightarrow{\simeq} & P \otimes_A N & \xrightarrow{\simeq} & Q \otimes_A N & \xrightarrow{\simeq} & \mathrm{N}(\mathcal{Q}\bar{\otimes}_\mathcal{A}\mathcal{N}) \\
\simeq \downarrow & & \downarrow & & \downarrow & & \downarrow \simeq \\
M' \otimes_A N & \longrightarrow & M \otimes_A N & =\!=\!= & M \otimes_A N & \longrightarrow & \mathrm{N}(\mathcal{M}\bar{\otimes}_\mathcal{A}\mathcal{N})\,.
\end{array}$$

The composition of the maps in the bottom line is the desired quasi-isomorphism. □

We interpolate a result from an earlier version of this paper, that is used in [17].

(2.3) *Künneth spectral sequences.* Let $\mathcal{A}$ be a simplicial ring, $\mathcal{M}$ a simplicial right $\mathcal{A}$-module, $\mathcal{N}$ a simplicial left $\mathcal{A}$-module, and let $A$, $M$, $N$ be the respective normalizations.

In a simplicial context, Quillen [36, §II.6] exhibits four Künneth spectral sequences that converge to the homotopy of the derived tensor product $\mathcal{M}\bar{\otimes}_\mathcal{A}^\mathbf{L}\mathcal{N}$; in particular

$${}^2\mathrm{E}_{p,q}^{\mathrm{Q}} = \mathrm{Tor}_p^{\pi(\mathcal{A})}(\pi(\mathcal{M}), \pi(\mathcal{N}))_q \implies \pi_{p+q}(\mathcal{M}\bar{\otimes}_\mathcal{A}^\mathbf{L}\mathcal{N})\,.$$

In a DG context, Eilenberg and Moore [34] construct a DG torsion product $\mathrm{Tor}^A(M, N)$ and approximate it by two spectral sequence, one of which has

$$^2\mathrm{E}^{\mathrm{EM}}_{p,q} = \mathrm{Tor}^{\mathrm{H}(A)}_p(\mathrm{H}(M), \mathrm{H}(N))_q \implies \mathrm{Tor}^A_{p+q}(M, N).$$

From the point of view of homotopical algebra [36], $\mathrm{Tor}^A(M, N) = \mathrm{H}(M \otimes^{\mathbf{L}}_A N)$, where $-\otimes^{\mathbf{L}}_A-$ is the derived tensor product on the category of DG modules over $A$.

By definition, the Quillen and Eilenberg-Moore spectral sequences have the same $^2\mathrm{E}$ page. The next statement was established in the course of the preceding proof.

PROPOSITION. *There is an isomorphism of spectral sequences $^r\omega\colon {^r\mathrm{E}^{\mathrm{EM}}} \longrightarrow {^r\mathrm{E}^{\mathrm{Q}}}$ with $^2\omega = \mathrm{id}$, that converges to an isomorphism of graded modules $\mathrm{H}(M \otimes^{\mathbf{L}}_A N) \cong \pi(\mathcal{M} \bar{\otimes}^{\mathbf{L}}_{\mathcal{A}} \mathcal{N})$.*

## 3. Deviations of local homomorphisms

In this section $\varphi\colon (R, \mathfrak{m}, k) \longrightarrow (S, \mathfrak{n}, \ell)$ is a local homomorphism.

Let $R \xrightarrow{\dot\varphi} R' \xrightarrow{\varphi'} \widehat{S}$ be a Cohen factorization of $\dot\varphi$. We denote $R'[Y]$ a DG algebra over $R'$ such that $R'[Y]^{\natural}$ is a tensor product of symmetric algebras of free modules with bases $Y_n$ for even $n \geq 0$ and exterior algebras of free modules with bases $Y_n$ for odd $n \geq 1$ is. Such a DG algebra is a *minimal model* of $\widehat{S}$ over $R'$ if $\mathrm{H}(R'[Y]) \cong \widehat{S}$, $Y = Y_{\geqslant 1}$, and the differential is decomposable in the sense that $\partial(Y) \subseteq \mathfrak{m}'R'[Y] + (Y)^2 R'[Y]$. Minimal models are characterized by the following properties: $Y = Y_{\geqslant 1}$; $\partial(Y_1)$ minimally generates the ideal $\mathfrak{a} = \mathrm{Ker}\,\varphi'$; $\{\mathrm{cls}(\partial(y)) \mid y \in Y_n\}$ minimally generates the $R'$-module $\mathrm{H}_{n-1}(R'[Y_{<n}])$ for $n \geq 2$; as a consequence, minimal models always exist, and have $Y_n$ finite for each $n$; for details we refer to [43] or [10, §7.2].

The next result shows that in the derived category of the category of $R$-algebras the isomorphism class of a minimal model is an invariant of $\dot\varphi$, and hence of the map $\varphi$:

(3.1) PROPOSITION. *If $R'[Y']$ and $R''[Y'']$ are minimal models of $\widehat{S}$ coming from Cohen factorizations of $\dot\varphi$, then there exist a minimal model $T[U]$ of $\widehat{S}$ coming from a Cohen factorization of $\dot\varphi$ and surjective quasi-isomorphisms*

$$R'[Y'] \longleftarrow T[U] \longrightarrow R''[Y'']$$

*of DG algebras over $R$ that induce the identity on $\widehat{S}$. Furthermore,*

$$\mathrm{card}\,(Y'_1) - \mathrm{edim}\,R' = \mathrm{card}\,(Y''_1) - \mathrm{edim}\,R''$$

*and*
$$\operatorname{card}(Y'_n) = \operatorname{card}(Y''_n) \quad \text{for } n \geq 2.$$

*Proof.* By [14, (1.2)] there exists a commutative diagram of ring homomorphisms

$$\begin{array}{ccccc}
 & & R' & & \\
 & \nearrow & \uparrow & \searrow & \\
R & \to & T & \to & \widehat{S}, \\
 & \searrow & \downarrow & \nearrow & \\
 & & R'' & &
\end{array}$$

where the horizontal row is a Cohen factorization and the vertical maps are surjections with kernels generated by $T$-regular sequences that extend to minimal sets of generators of the maximal ideal of $T$. Thus, we may assume that there is a surjective homomorphism $R'' \to R'$ with kernel of this type, and switch the notation accordingly.

Changing $Y''_1$ if necessary, we may also assume that it contains a subset $V$ such that $\partial(V)$ minimally generates $\operatorname{Ker}(R'' \to R')$. As $\partial(V)$ is a regular sequence, the Koszul complex $R''[V]$ has $\operatorname{H}(K) \cong R'$, and is a DG subalgebra of $R''[Y'']$. Since $R''[Y'']^\natural$ is a free module over $R''[V]^\natural$, we conclude by (2.1.1) that the canonical map $R''[Y''] \to R''[Y'']/(\partial(V), V) = R'[\overline{Y'}]$, where $\overline{Y'} = Y'' \smallsetminus V$, is a quasi-isomorphism. Thus, $\operatorname{H}(R'[\overline{Y'}]) \cong \widehat{S}$.

The differential of $R'[\overline{Y'}]$ inherits the decomposability of that of $R''[Y'']$, so $R'[\overline{Y'}]$ is a minimal model of $\widehat{S}$ over $R'$. By [10, (7.2.3)] the DG algebras $R'[\overline{Y'}]$ and $R'[Y]$ are isomorphic over $R'$; hence $\overline{Y'}_n = Y_n$ for all $n$. Now note that $\operatorname{card} Y'_1 = \operatorname{card} \overline{Y'}_1 = \operatorname{card}(Y''_1) - (\operatorname{edim} R'' - \operatorname{edim} R')$, and $\operatorname{card}(Y'_n) = \operatorname{card} \overline{Y'}_n = \operatorname{card}(Y''_n)$ for $n \geq 2$. □

In view of the proposition we refer to a minimal model of $\widehat{S}$ over the ring $R'$ in any Cohen factorization of $\dot{\varphi}$ as a *minimal model of the homomorphism* $\dot{\varphi}$. We call the number

$$\varepsilon_n(\varphi) = \begin{cases} \operatorname{card}(Y_1) - \operatorname{edim} R' + \operatorname{edim} S/\mathfrak{m}S & \text{for } n = 2; \\ \operatorname{card}(Y_{n-1}) & \text{for } n \geq 3, \end{cases}$$

the $n^{\text{th}}$ *deviation* of $\varphi$. To explain the terminology, note that if the ring $R$ is regular and $\varphi$ is surjective then [10, (7.2.7)] shows that $\varepsilon_n(\varphi) = \varepsilon_n(S)$ for $n \geq 2$, where the $n^{\text{th}}$ deviation $\varepsilon_n(S)$ of the local ring $S$ is classically defined in terms of an infinite product decomposition of its Poincaré series $\sum_{n=0}^{\infty} \operatorname{rank}_\ell \operatorname{Tor}_n^S(\ell, \ell)$.

The deviations of a local ring measure its failure to be regular, or c.i. The vanishing of the initial deviations of a local homomorphisms are interpreted along similar lines.

(3.2) *Remark.* An equality $\varepsilon_2(\varphi) = 0$ means that $\varphi$ is flat with $S/\mathfrak{m}S$ regular, and is equivalent to the vanishing of $\varepsilon_n(\varphi)$ for $n \geq 2$.

Indeed, by [14, (1.5)] there is a Cohen factorization with $\operatorname{edim} R' = \operatorname{edim} S$. If $\varepsilon_2(\varphi) = 0$ then $Y_1 = \varnothing$, so $\widehat{S} = \operatorname{H}_0(R'[Y]) = R'$, hence $\widehat{S}$ is flat over $R$ and $\widehat{S}/\mathfrak{m}\widehat{S}$ is regular; these properties descend to $S$ and $S/\mathfrak{m}S$. Conversely, if $\varphi$ is flat with regular closed fiber, then $R \to \widehat{S} = \widehat{S}$ is a Cohen factorization of $\grave{\varphi}$, so $\grave{\varphi}$ has a minimal model with $Y = \emptyset$.

(3.3) *Remark.* An equality $\varepsilon_3(\varphi) = 0$ means that $\varphi$ is c.i. at $\mathfrak{n}$, and is equivalent to the vanishing of $\varepsilon_n(\varphi)$ for $n \geq 3$; as a consequence, if $\varphi$ is c.i. at $\mathfrak{n}$ then in *each* Cohen factorization of $\grave{\varphi}$ the kernel of the surjective map $\varphi'$ is generated by a regular sequence.

Indeed, the definitions of c.i. homomorphism and of deviations of a homomorphism allow us to replace $\varphi$ by $\varphi'$; changing notation, we may assume that $\varphi \colon R \to S$ is surjective. In this situation $\varepsilon_3(\varphi) = \operatorname{card}(Y_2)$ is the minimal number of generators of $\operatorname{H}_1(R[Y_1])$, where $R[Y_1]$ is the Koszul complex on a minimal set of generators of $\mathfrak{a} = \operatorname{Ker} \varphi$. Thus $\varepsilon_3(\varphi)$ vanishes if and only if $\mathfrak{a}$ is generated by a regular sequence, that is, if and only if $\varphi$ is c.i. at $\mathfrak{n}$. When this is the case the Koszul complex is exact, so $R[Y] = R[Y_1]$; in other words, we have $\varepsilon_n(\varphi) = \operatorname{card}(Y_{n-1}) = 0$ for $n \geq 3$.

We establish the *rigidity* of deviations for homomorphisms of finite flat dimension, strengthening a result of Avramov and Halperin [16]: If $\varphi$ is not c.i. at $\mathfrak{n}$, then $\varepsilon_n(\varphi) \neq 0$ for $n \gg 0$ (it is stated there for 'factorizable' homomorphisms, but the construction of Cohen factorizations in [14] shows that each $\varphi$ has this property).

(3.4) THEOREM. *If* $\operatorname{fd}_R S < \infty$ *and* $\varepsilon_n(\varphi) = 0$ *for some* $n \geq 4$, *then* $\varphi$ *is c.i. at* $\mathfrak{n}$.

When $R$ is regular and $\varphi$ is surjective, the theorem is equivalent to Halperin's result [30] on the rigidity of deviations of local rings. We extend his argument to the relative case by using Cohen factorizations, and develop shortcuts based on the study of derivations in [10]. First, we record how conditions on the flat dimension of $\varphi$ pass through factorizations.

(3.5) *Remark.* If $R \to R' \to \widehat{S}$ is a Cohen factorization of $\grave{\varphi}$ then

$$\operatorname{fd}_R S = \operatorname{fd}_R \widehat{S} \leq \operatorname{pd}_{R'} \widehat{S} \leq \operatorname{fd}_R S + \operatorname{edim}(S/\mathfrak{m}S)$$

where edim $R$ denotes the minimal number of generators of $\mathfrak{m}$ and $\mathrm{pd}_{R'}\widehat{S}$ is the projective dimension of the $R'$-module $\widehat{S}$ (cf. [14, (3.2)] or [16, (3.2)]).

In particular, $\mathrm{fd}_R S$ is finite only if $\mathrm{pd}_{R'}\widehat{S}$ is finite.

We recall some basics on Tate's [41] construction of DG algebra resolutions (for details, see [29], [10]). When $A$ is a DG algebra $A\langle X\rangle$ denotes a DG algebra obtained from it by adjunctions of sets of exterior variables $X_n$ in odd degrees $n \geq 1$ and of divided power variables in even degrees $n \geq 2$. The $i^{\mathrm{th}}$ divided power of $x \in X_{\mathrm{even}}$ is denoted $x^{(i)}$. It satisfies, among other relations, $|x^{(i)}| = i|x|$; $x^{(0)} = 1$; $x^{(1)} = x$, as well as

$$x^{(i)}x^{(j)} = \binom{i+j}{i}x^{(i+j)} \qquad \text{and} \qquad \partial(x^{(i)}) = \partial(x)x^{(i-1)} \qquad \text{for all} \quad i, j \geq 0.$$

We say that $X$ is a set of $\Gamma$-*variables* adjoined to $A$ and $A\langle X\rangle$ is a $\Gamma$-*free extension* of $A$.

(3.6) *Remark.* If $A_0$ is a local ring with maximal ideal $\mathfrak{m}$ and residue field $\ell$, then $A\langle X\rangle$ is an *acyclic closure* of $\ell$ over $A$ if $X = X_{\geqslant 1}$ and $\partial$ satisfies the conditions: $\partial(X_1)$ minimally generates $\mathfrak{m}$ modulo $\partial(A_1)$ and the classes of $\{\partial(x) \mid x \in X_n\}$ minimally generate the $A_0$-module $\mathrm{H}_{n-1}(A\langle X_{<n}\rangle)$ for $n \geq 2$. Gulliksen [29, (6.2))] proves that if $A\langle X\rangle$ is an acyclic closure of $\ell$, then $\partial(A\langle X\rangle) \subseteq (\mathfrak{m} + A_{\geqslant 1})A\langle X\rangle$ (cf. also [10, (6.3.4)]).

We need a simple case of [10, (7.2.11)].

(3.7) LEMMA. *A DG algebra $\ell[Y]$ with $Y = Y_{\geqslant 1}$ and $\partial(Y) \subseteq (Y)^2\ell[Y]$ has a $\Gamma$-free extension $B = \ell[Y]\langle X\rangle$ with $X = \{x_y \mid y \in Y, |x_y| = |y| + 1\}$, $\mathrm{H}(B) = \ell$, and $\partial(B) \subseteq (Y)B$.*

*Proof.* Set $\ell[Y_{\geqslant n}] = \ell[Y]/(Y_{<n})$. Starting with $B^0 = \ell[Y]$ and $X_{\leqslant 0} = \emptyset$, assume by induction that for some $n \geq 0$ we have a surjective quasi-isomorphism $\varkappa^n \colon B^n = \ell[Y]\langle X_{\leqslant n}\rangle \longrightarrow \ell[Y_{\geqslant n}]$, where $X_i = \{x_y \mid y \in Y_{i-1}\}$ denotes a set of $\Gamma$-variables of degree $i$. The condition $\partial(Y) \subseteq (Y)^2\ell[Y]$ implies that $Y_n$ is a basis of $\mathrm{H}_n(\ell[Y_{\geqslant n}])$ over $\ell$. Thus, for each $y \in Y_n$ there is a cycle $z_y \in \mathrm{Z}_n(\ell[Y]\langle X_{\leqslant n}\rangle)$ such that $\varkappa^n(z_y) = y$. Choosing a set $X_{n+1} = \{x_y \mid y \in Y_n\}$ of $\Gamma$-variables of degree $n + 1$, we extend $\varkappa^n$ to a morphism

$$B^{n+1} = \ell[Y]\langle X_{\leqslant n}\rangle\langle X_{n+1} \mid \partial(x_y) = z_y\rangle \longrightarrow \ell[Y_{\leqslant n}]\langle X_{n+1} \mid \partial(x_y) = y\rangle = C^{n+1}$$

of DG algebras that is the identity on $X_{n+1}$; it is easily seen to be a quasi-isomorphism. It is well known that the DG subalgebra

$$D^{n+1} = \ell[Y_n]\langle X_{n+1} \mid \partial(x_y) = y\rangle$$

of $C^{n+1}$ has $\mathrm{H}(D^{n+1}) = \ell$ (cf. Cartan [18]). From (2.1.1) we see that $C^{n+1} \longrightarrow C^{n+1} \otimes_{D^{n+1}} \ell = \ell[Y_{\geqslant n+1}]$ is a quasi-isomorphism. In the limit we get a quasi-isomorphism $\mathrm{inj}\lim_n \varkappa^n \colon \mathrm{inj}\lim_n B^n \to \mathrm{inj}\lim_n \ell[Y_{\leqslant n}]$, which is just the canonical augmentation $B \to \ell$. Since we have constructed $B$ as an acyclic closure of $\ell$ over $\ell[Y]$, we have $\partial(B) \subseteq (Y)B$ by (3.6). □

(3.8) LEMMA. *Let* $R \to R' \to \widehat{S}$ *be a Cohen factorization of* $\dot{\varphi}$, *let* $R'[Y]$ *be a minimal model of* $\widehat{S}$, *and set* $\ell[Y_{\geqslant n}] = R'[Y]/(\mathfrak{m}', Y_{<n})R'[Y]$ *for* $n \geq 1$.

*If* $\mathrm{fd}_R S < \infty$, *then for each* $n \geq 1$ *the product of any* $q$ *elements of positive degree in* $\mathrm{H}(\ell[Y_{\geqslant n}])$ *is trivial when* $q = \mathrm{fd}_R S + \mathrm{edim}\,(S/\mathfrak{m}S) + 1$.

*Proof.* The DG algebra $\ell[Y] = R'[Y] \otimes_{R'} \ell$ has $\mathrm{H}_i(\ell[Y]) \cong \mathrm{Tor}_i^{R'}(S,\ell) = 0$, and by (3.5) this module is trivial when $i \geq q$. Setting $J_i = 0$ for $i \leq q-2$, $J_{q-1} = \partial(\ell[Y]_q)$, and $J_i = \ell[Y]_i$ for $i \geq q$, we get a subcomplex $J \subseteq \ell[Y]$ with $\mathrm{H}(J) = 0$; for degree reasons, it is a DG ideal of $\ell[Y]$, so $\ell[Y] \to C = \ell[Y]/J$ is a quasi-isomorphism of DG algebras.

Let $B = \ell[Y]\langle X \rangle$ be the $\Gamma$-free extension of Lemma (3.7). We set $B^n = \ell[Y]\langle X_{\leqslant n}\rangle$ and prove that $\bigl(\mathrm{H}_{\geqslant 1}(B^n)\bigr)^q = 0$. If $n = 0$ then $B^0 = \ell[Y]$ is exact in degrees $\geq q$ and the assertion is clear. If $n > 0$, then due to $\mathrm{H}(B) = \ell$ and $\partial(B) \subseteq (Y)B$ we have

$$\mathrm{Z}_{\geqslant 1}(B^n) = B^n \cap \mathrm{Z}_{\geqslant 1}(B) = B^n \cap \partial(B) \subseteq B^n \cap (Y)B = (Y)B^n.$$

As $B^{n\natural}$ is a free module over $\ell[Y]^\natural$ the canonical map $B^n \to B^n/JB^n$ is a quasi-isomorphism by (2.1.1), so $\mathrm{H}(JB^n) = 0$. In view of the preceding computation, this implies

$$(\mathrm{Z}_{\geqslant 1}(B^n))^q \subseteq \mathrm{Z}(B^n) \cap (Y)^q B^n \subseteq \mathrm{Z}(B^n) \cap JB^n = \mathrm{Z}(JB^n) = \partial(JB^n).$$

We conclude that $\bigl(\mathrm{H}_{\geqslant 1}(B^n)\bigr)^q = 0$ and finish the argument by invoking the quasi-isomorphism $B^n = \ell[Y]\langle X_{\leqslant n}\rangle \to \ell[Y_{\geqslant n}]$ established in the preceding proof. □

Let $A\langle X \rangle$ be an extension of a DG algebra $A$ by a set of $\Gamma$-variables $X = X_{\geqslant 1}$, and let $U$ be a DG module over $A\langle X \rangle$. A (chain) $A$-*linear $\Gamma$-derivation* is a homogeneous (chain) map $\vartheta \colon A\langle X \rangle \to U$, such that the relations

$$\vartheta(a) = 0, \qquad \vartheta(bb') = \vartheta(b)b' + (-1)^{|b||b'|}\vartheta(b')b, \qquad \vartheta(x^{(i)}) = x^{(i-1)}\vartheta(x)$$

hold for all $a \in A$, $b, b' \in A\langle X \rangle$, $x \in X_{\mathrm{even}}$, and $i \in \mathbb{N}$.

Let $\mathrm{H}_0(A\langle X \rangle) = S$ and set $\mathfrak{a} = \mathrm{Ker}\,(A_0 \to S)$. If $X^{(2)}$ denotes the set of all products $x_r^{(i_r)} \cdots x_s^{(i_s)}$ with $i_r + \cdots + i_s \geq 2$ then $D = A + \mathfrak{a}X + AX^{(2)}$ is a DG submodule of $A\langle X \rangle$, so the canonical projection $\pi \colon A\langle X \rangle \to L = A/D$ makes $L$ into a complex of free $S$-modules, with $X_n$ a basis of $L_n$ for each $n$. We call $L$ the *complex of indecomposables* of the extension $A \to A\langle X \rangle$. The following is proved in [10, (6.3.6)].

(3.9) *Remark.* Let $U$ be a DG module over $A$ with $U_i = 0$ for $i < 0$ and let $\beta\colon U \to M$ be a surjective quasi-isomorphism to a complex of $S$-modules.

Each chain map $\xi\colon L \to M$ of degree $-n$ lifts to a chain $\Gamma$-derivation $\vartheta\colon A\langle X\rangle \to U$ of degree $-n$, such that $\beta\vartheta = \xi\pi$; if a family $\{u_x \in U_0\}_{x \in X_n}$ satisfies $\beta(u_x) = \xi(x)$ for $u_x \in X_n$, then $\xi$ may be chosen with $\xi(x) = u_x$ for each $x \in X_n$.

*Proof of Theorem* (3.4). Assume that there exists a local homomorphism $\varphi$ such that $\operatorname{fd}_R S < \infty$, and $\varepsilon_n(\varphi) = 0$ for some $n \geq 4$, but which is not c.i. at $\mathfrak{n}$. By Remark (3.5) we may further assume that $\varphi$ is surjective and $\operatorname{pd}_R S < \infty$. Fix a minimal model $R[Y]$ of $\varphi$ and for each $n \geq 1$ set $\ell[Y_{\geq n}] = R[Y]/(\mathfrak{m}, Y_{<n})$, where $\ell$ is the residue field of $R$.

Since $\varepsilon_3(\varphi) \neq 0$ by Remark (3.3), we can find $j$ with $\varepsilon_j(\varphi) \neq 0 = \varepsilon_{j+1}(\varphi)$, so that $Y_{j-1} \neq \emptyset$ and $Y_j = \emptyset$. Let $i$ be the integer part of $\frac{j-1}{2}$, set $A = R[Y_{\leq 2i-1}]$, and let $A\langle X\rangle$ be an acyclic closure of $S = \operatorname{H}_0(A)$ over $A$. Here is the key

*Claim.* There exists an $\ell$-linear chain $\Gamma$-derivation $\theta$ of $\ell\langle X\rangle = \ell \otimes_A A\langle X\rangle$, such that $\theta(x) = 1$ for some cycle $x \in X_{2i}$.

It implies that $\operatorname{cls}(x^{(r)}) \neq 0 \in \operatorname{H}(\ell\langle X\rangle)$ for all $r \geq 0$. Indeed, if $x^{(r)} = \partial(v)$, then $1 = \theta^r(x^{(r)}) = \theta^r\partial(v) = \partial\theta^r(v) = 0$, which is absurd. The multiplication table for divided powers then shows that $x^r = r!x^{(r)} \neq 0$ when $\operatorname{char}(\ell) = 0$, and that $x \cdot x^{(p)} \cdots x^{(p^r)} = x^{(1+p+\cdots+p^r)} \neq 0$ when $\operatorname{char}(\ell) = p > 0$.

On the other hand, the graded algebra underlying $R[Y] = A[Y_{\geq 2i}]$ is a free extension of $A^\natural$, and $A\langle X\rangle \to S$ is a surjective quasi-isomorphism, so the inclusion $R[Y] = A \subseteq A\langle X\rangle$ extends to a morphism of DG algebras $\phi\colon R[Y] \to A\langle X\rangle$; it is necessarily a quasi-isomorphism. By (2.1.2), so is the induced map

$$\ell[Y_{\geq 2i}] = \ell \otimes_A A[Y_{\geq 2i}] \xrightarrow{\ell \otimes_A \phi} \ell \otimes_A A\langle X\rangle = \ell\langle X\rangle\,.$$

By Lemma (3.8) the product of every $(\operatorname{pd}_R S + 1)$ elements of positive degree in $\operatorname{H}(\ell\langle X\rangle)$ is trivial: we have a contradiction, so it remains to establish the claim.

For this, we take a closer look at the construction of $\phi\colon A[Y_{\geq 2i}] \to A\langle X\rangle$.

Since $\operatorname{H}_{n-1}(A) = 0$ for $1 < n < 2i$, we have $X_n = \emptyset$ for $n < 2i$, and we can take $X_{2i} = \{x_y \mid y \in Y_{2i}\}$, with $\partial(x_y) = \partial(y)$. Let $\phi^{2i}\colon A[Y_{2i}] \to A\langle X_{2i}\rangle$ be the morphism of DG algebras over $A$, such that $\phi^{2i}(y) = x_y$ for $y \in Y_{2i}$. The map $\phi_n^{2i}$ is bijective for $n < 4i$; hence $\operatorname{H}_{2i}(\phi^{2i})$ is an isomorphism. By construction, the classes of the cycles $\partial(Y_{2i+1})$ minimally generate $\operatorname{H}_{2i}(A[Y_{2i}])$.

When $j = 2i + 1$, we have $\operatorname{H}_{2i}(A[Y_{2i}]) = 0$, and hence $X_{2i+1} = \emptyset$. When $j = 2i + 2$, we take $X_{2i+1} = \{x_y \mid y \in Y_{2i+1}\}$ and extend the map $y \mapsto x_y$ to a morphism $\phi^{2i+1}\colon R[Y_{\leq 2i+1}] = A[Y_{2i}][Y_{2i+1}] \to A\langle X_{\leq 2i+1}\rangle$ of DG algebras over $A[Y_{2i}]$; for the same reasons as above, we conclude that $\phi_n^{2i+1}$ is bijective for $n < 4i$. Thus, if $i \geq 2$ then $0 = \mathrm{H}_{2i+1}(R[Y_{\leqslant 2i+1}]) \cong \mathrm{H}_{2i+1}(A\langle X_{\leqslant 2i+1}\rangle)$, and so $X_{2i+2} = \emptyset$; if $i = 1$ (that is, if $j = 4$), then $\phi^3$ still induces an isomorphism $\mathrm{Z}_3(R[Y_{\leqslant 3}]) \cong \mathrm{Z}_3(A\langle X_{\leqslant 3}\rangle)$, hence a surjection $0 = \mathrm{H}_3(R[Y_{\leqslant 3}]) \to \mathrm{H}_3(A\langle X_{\leqslant 3}\rangle)$ that gives $X_4 = \emptyset$.

As a result of the preceding discussion we now know that

$$X_n = \emptyset \quad \text{for} \quad n < 2i\,; \qquad X_{2i} \neq \emptyset\,; \qquad \begin{aligned} X_{2i+1} &= \emptyset \quad \text{if} \quad j = 2i+1\,; \\ X_{2i+2} &= \emptyset \quad \text{if} \quad j = 2i+2\,. \end{aligned}$$

We arbitrarily pick $x \in X_{2i}$, and set $X'_{2i} = X_{2i} \smallsetminus \{x\}$.

When $j = 2i + 1$ the complex of indecomposables of Remark (3.9) has the form

$$L\colon \ldots \to SX_{2i+2} \to 0 \to Sx \oplus SX'_{2i} \to 0 \to \ldots$$

so $\mathrm{H}_{2i}(L) = SX_{2i}$ is a free $S$-module. If $M$ is the complex of $S$-modules with $S$ concentrated in degree 0, then $\xi(x) = 1$ and $\xi(X \smallsetminus \{x\}) = 0$ defines a chain map $\xi\colon L \to M$ of degree $-2i$. By Remark (3.9) it lifts over the augmentation $\varepsilon\colon A\langle X\rangle \to S$ to an $A$-linear chain $\Gamma$-derivation $\vartheta\colon A\langle X\rangle \to A\langle X\rangle$ with $\vartheta(x) = 1$, so $\theta = \ell \otimes_A \vartheta$ acts as desired.

When $j = 2i + 2$ the complex of indecomposables has the form

$$L\colon \ldots \to SX_{2i+3} \to 0 \to SX_{2i+1} \xrightarrow{\partial_{2i+1}} Sx \oplus SX'_{2i} \to 0 \to \ldots\,.$$

Consider the complex of free $S$-modules, concentrated in degrees 0 and 1

$$M\colon \ldots \to 0 \to S^e \xrightarrow{\delta} S \to 0 \to \ldots$$

with $\delta = (t_1, \ldots, t_e)$, where $\{t_1, \ldots, t_e\}$ is a minimal set of generators of $\mathfrak{n}$. Let $\xi_{2i}\colon L_{2i} \to M_0$ be the $S$-linear map defined by $\xi_{2i}(x) = 1$ and $\xi_0(X'_{2i}) = 0$. Since $\mathrm{Im}(\partial_{2i+1}) \subseteq \mathfrak{n}X_{2i}$ by Remark (3.6), we have $\mathrm{Im}(\xi_{2i}\partial_{2i+1}) \subseteq \mathfrak{n}M_0 = \mathrm{Im}\,\delta$, so there is an $S$-linear homomorphism $\xi_{2i+1}\colon L_{2i+1} \to M_1$ with $\delta\xi_{2i+1} = \xi_{2i}\partial_{2i+1}$. Setting $\xi_n = 0$ for $n \geq 2i + 2$, we get a chain map $\xi\colon L \to M$ of degree $-2i$.

Let $U$ be the DG module over $A\langle X\rangle$ with $U^\natural = (A\langle X\rangle)^\natural \oplus \bigoplus_{h=1}^e (A\langle X\rangle u_h)^\natural$ and $\partial(u_h) = t_h$ for $h = 1, \ldots, e$. The augmentation $\varepsilon\colon A\langle X\rangle \to S$ induces a quasi-isomorphism $\beta\colon U \to U \otimes_{A\langle X\rangle} S = M$, so Remark (3.9) yields an $A$-linear $\Gamma$-derivation $\vartheta\colon A\langle X\rangle \to U$ with $\vartheta(x) = 1$. Now $\overline{U} = \ell \otimes_A U = \ell\langle X\rangle \oplus \bigoplus_{h=1}^e \ell\langle X\rangle u_h$ is a DG module over $\ell\langle X\rangle$ with $\partial(u_h) = 0$ for $h = 1, \ldots, e$. The composition $\theta$ of $\ell \otimes_A \vartheta\colon \ell\langle X\rangle \to \overline{U}$ with the projection $\overline{U} \to \overline{U}/(u_1, \ldots, u_e) = \ell\langle X\rangle$ has the desired properties. □

Finally, we establish the *exponential growth* of deviations: this is a 'looking glass' version, in the sense of [8] and [15], of a result of Félix, Halperin, and Thomas [23] on the ranks of homotopy groups of finite CW complexes. The strategy of their proof goes through, but its engine—the 'mapping theorem' of Félix and Halperin [22] for rational Ljusternik-Schnirelmann category—is not available in our context.

(3.10) THEOREM. *If* $\mathrm{fd}_R S < \infty$ *and* $\varphi$ *is not c.i. at* $\mathfrak{n}$, *then there exist a real number* $\gamma > 1$ *and a sequence of positive integers* $s_j$ *with* $j \geq 0$, *such that* $2s_j \leq s_{j+1} \leq qs_j$ *for* $q = \mathrm{fd}_R S + \mathrm{edim}(S/\mathfrak{m}S) + 1$ *and*

$$\varepsilon_{s_j}(\varphi) > \gamma^{s_j} \qquad \text{for} \quad j \geq 1.$$

*Proof.* As in the preceding proof, we first put ourselves in a situation when $\varphi = \grave{\varphi}$ is surjective (so $k = \ell$), $\mathrm{pd}_R S$ is finite, $R[Y]$ is a minimal model of $\varphi$, and $\ell[Y_{\geqslant n}] = R[Y]/(\mathfrak{m}, Y_{<n})$ for $n \geq 1$. We write $Y_{[n]}$ for the span of $\bigcup_{j=n}^{2n} Y_j$. Letting $Y_n$ denote also the $\ell$-linear span of the variables $y \in Y_n$, we can write $Y_{[n]}^i$ for the span of all products involving $i$ elements of $Y_{[n]}$. Finally, we set $a_n = \mathrm{card}\,(Y_n)$ and $s(n) = \sum_{j=n}^{2n} \mathrm{card}\,(Y_j)$.

The first goal is to prove that the sequence $\{a_n\}_{n \geqslant 2}$ is unbounded.

For every $y \in Y_{\geqslant n}$ there are uniquely defined $\alpha_i(y) \in Y_{[n]}^i \subseteq \ell[Y_{\geqslant n}]$, such that

$$\partial(y) \equiv \sum_{i \geqslant 2} \alpha_i(y) \mod \left((Y_{>2n})\ell[Y_{\geqslant n}]\right).$$

For each $i$, the assignment $y \mapsto \alpha_i(y)$ defines an $\ell$-linear homomorphism $\alpha_i \colon Y_{\geqslant n} \to Y_{[n]}^i$. By the decomposability of the differential of $\ell[Y_{\geqslant n}]$ we have $\partial(Y_{[n]}^i) = 0$ for all $i \geq 2$. For $q = \mathrm{fd}_R S + \mathrm{edim}\,(S/\mathfrak{m}S) + 1$, Lemma (3.8) shows that $Y_{[n]}^q$ consists of boundaries; hence

(3.10.1)
$$\sum_{i=2}^{q} Y_{[n]}^{q-i} \alpha_i(Y_{\geqslant n}) \supseteq Y_{[n]}^q.$$

For degree reasons, $\alpha_i(Y_j) = 0$ when $j < in + 1$ or $j > i(2n) + 1$, so

(3.10.2) $\quad s(in+1) = \displaystyle\sum_{j=in+1}^{2in+2} a_j \geq \sum_{j=in+1}^{i(2n)+1} \mathrm{rank}_\ell\, \alpha_i(Y_j) = \mathrm{rank}_\ell\, \alpha_i(Y_{\geqslant n})\,.$

Set $d = (2q)^q$ and $n_0 = 2(qd)^2$. By Theorem (3.4), we have $s(n_0) > (qd)^2$, so assume by induction that integers $n_0, n_1, \ldots, n_j$ have been found with
(3.10.3)
$\quad q(n_{h-1} + 1) \geq n_h + 1 \qquad \text{and} \qquad s(n_h) \geq (qd)s(n_{h-1}) \quad \text{for } 1 \leq h \leq j\,.$

Choose $n_{j+1} = in_j + 1$ such that $s(n_{j+1}) = \max\{s(in_j + 1) \mid 2 \leq i \leq q\}$. It is then clear that $q(n_j + 1) \geq n_{j+1} + 1$. Using (3.10.2) and (3.10.1), we get

$$(q-1)s(n_j)^{q-2}s(n_{j+1}) \geq \sum_{i=2}^{q} s(n_j)^{q-i} s(in_j + 1)$$

$$\geq \sum_{i=2}^{q} \left(\operatorname{rank}_\ell Y_{[n_j]}\right)^{q-i} \operatorname{rank}_\ell \alpha_i(Y_{\geqslant n_j})$$

$$\geq \operatorname{rank}_\ell \left(\sum_{i=2}^{q} Y_{[n_j]}^{q-i} \alpha_i(Y_{\geqslant n_j})\right) \geq \operatorname{rank}_\ell Y_{[n_j]}^q \geq \binom{s(n_j)}{q}$$

$$\geq \frac{s(n_j)^q}{(2q)^q} = \frac{s(n_j)}{d} s(n_j)^{q-1} \geq (q^2 d) s(n_j)^{q-1}$$

so $s(n_{j+1}) \geq (qd)s(n_j)$, completing the induction step.

From (3.10.3) we see that $q^j(n_0 + 1) \geq n_j + 1$ and $s(n_j) \geq q^j s(n_0) d^j$ for $j \geq 1$. Assuming that $c \geq a_n$ for some $c$ and all $n$, we get $c(n_j + 1) \geq s(n_j)$, hence $c(n_0 + 1) \geq s(n_0) d^j$ for all $j \geq 1$. This is impossible, so the sequence $\{a_n\}_{n \geq 2}$ is unbounded, as desired.

Set $b = (2q)^{q+1}$, choose $r_1$ so that $a_{r_1} = a > b$, and assume by induction that $r_1, \ldots, r_j$ have been found, with $r_h = i_{h-1}r_{h-1} + 1$ for $2 \leq i_{h-1} \leq q$ and $a_{r_h} \geq a_{r_{h-1}}^{i_{h-1}}/b$ for $2 \leq h \leq j$. The condition $\beta(y) \equiv \partial(y) \mod \left((Y_{>r_j})\ell[Y_{\geqslant r_j}]\right)$ defines an $\ell$-linear map $\beta \colon Y_{\geqslant r_j} \to \sum_{i \geqslant 2} Y_{r_j}^i$. Since $\beta(y) = 0$ unless $|y| \equiv 1 \pmod{r_j}$, using the decomposability of the differential and Lemma (3.8) as before, we get $\sum_{i=2}^{q} Y_{r_j}^{q-i} \beta(Y_{ir_j+1}) \supseteq Y_{r_j}^q$. It follows that

$$\sum_{i=2}^{q} a_{r_j}^{q-i} a_{ir_j+1} \geq \binom{a_{r_j}}{q} \geq \frac{a_{r_j}^q}{(2q)^q} = (2q)\frac{a_{r_j}^q}{b} \,;$$

hence $a_{ir_j+1} \geq a_{r_j}^i/b$ for some $i$ with $2 \leq i \leq q$; set $i_j = i$ and $r_{j+1} = i_j r_j + 1$.

The induction is now complete, so for each $j \geq 1$ we have an expression $r_{j+1} = u_j r_1 + v_j$, with $v_j = (i_j i_{j-1} \cdots i_2) + \cdots + (i_j i_{j-1}) + i_j + 1$ and $u_j = i_j i_{j-1} \cdots i_1$. We note that

$$u_j > \left(\frac{1}{2} + \cdots + \frac{1}{2^j}\right) u_j \geq \left(\frac{1}{i_1} + \cdots + \frac{1}{i_1 \cdots i_j}\right) u_j = v_j,$$

and hence $u_j(r_1 + 1) > u_j r_1 + v_j = r_{j+1}$. Since $a > b > 1$, we obtain

$$a_{r_{j+1}} \geq \frac{a_{r_j}^{i_j}}{b} \geq \frac{a_{r_{j-1}}^{i_j i_{j-1}}}{b^{1+i_j}} \geq \cdots \geq \frac{a^{u_j}}{b^{v_j}} > \left(\frac{a}{b}\right)^{u_j} > \left(\frac{a}{b}\right)^{\frac{r_{j+1}}{r_1+1}}.$$

Setting $s_j = r_j + 1$ and $\gamma = \sqrt[s]{a/b}$, we see that $\gamma > 1$ and that

$$\varepsilon_{s_j}(\varphi) = \operatorname{card}(Y_{r_j}) = a_{r_j+1} > \gamma^{r_j+1} = \gamma^{s_j}$$
$$2s_j = 2r_j + 2 \leq r_{j+1} + 1 = s_{j+1} \leq (qr_j + 1) + 1 = qs_j - q + 2 \leq qs_j$$

for all integers $j \geq 1$. This is the desired conclusion. □

## 4. André-Quillen homology

Let $\varphi' \colon R \to S'$ and $S' \to k'$ be homomorphisms of rings, with $\varphi'$ surjective. Quillen links the cotangent invariants of $\varphi'$ to invariants defined in terms of classical derived functors by a fundamental first quadrant homological spectral sequence of $k'$-algebras

$$^2\mathrm{E}_{p,q} = \pi_{p+q}(\operatorname{Sym}_q^{k'}(\mathcal{L}(S'|R) \otimes_{S'} k')) \implies \operatorname{Tor}_{p+q}^R(S', k')$$

where $\operatorname{Sym}^{k'}$ denotes the symmetric algebra functor over $k'$, applied dimensionwise, and $\mathcal{L}(S'|R)$ is a simplicial $S'$-module that underlies the construction of the cotangent complex, in the sense that $\mathrm{L}(S|R) = \mathrm{N}(\mathcal{L}(S|R))$ (cf. [37, (6.3)]); note that the $^2\mathrm{E}$ page only depends on $\pi(\mathcal{L}(S|R) \otimes_S k') = \mathrm{D}(S|R, k')$ (cf. Dold [19]).

We are interested in an edge homomorphism, for which we review some more divided powers. A DG $\Gamma$-*algebra* is a graded strictly commutative algebra in which each element $x$ of even positive degree has a system of divided powers $x^{(i)}$ compatible with the differential (cf. [18], [29]). It is well known—but a complete reference is hard to find—that the normalization $\mathrm{N}\mathcal{A}$ of a commutative simplicial ring $\mathcal{A}$ is a DG $\Gamma$-algebra, and the divided powers pass to the homotopy $\pi(\mathcal{A}) = \mathrm{H}(\mathrm{N}\mathcal{A})$. (Cartan's text [18] is the original source, but the statement does not seem to appear there; $x^{(i)}$ can be defined directly in terms of Eilenberg-Zilber shuffles, as in Nicollerat [35, p. 660] or Goerss [25, p. 30].)

(4.1) *Remark.* If $\mathcal{S}$ is a cofibrant simplicial algebra resolution of the $R$-algebra $S'$ (cf. [37]), then $\operatorname{Tor}^R(S', k') = \mathrm{H}(\mathrm{N}(\mathcal{S} \bar{\otimes}_R k'))$ is a $\Gamma$-algebra by the preceding remarks. Let $\mathrm{T}_n(R, S', k')$ be the quotient of $\operatorname{Tor}^R(S', k')$ by its graded $k'$-submodule spanned by 1, all products $uv$ with $u, v \in \operatorname{Tor}_{\geqslant 1}^R(S', k')$, and all divided powers $w^{(i)}$ with $w \in \operatorname{Tor}_{2j}^R(S', k')$ for $i \geq 2$, and $j \geq 1$. By [37, (6.5)] or [2, §6], the edge homomorphism

$$\operatorname{Tor}_n^R(S', k') \twoheadrightarrow {}^\infty\mathrm{E}_{n-1,1} \hookrightarrow {}^2\mathrm{E}_{n-1,1} = \mathrm{D}_n(S'|R, k')$$

induces $k'$-linear maps $\beta_n \colon \mathrm{T}_n(R, S', k') \to \mathrm{D}_n(S'|R', k')$. When $S'$ and $k'$ are fields of characteristic 0 Quillen [37, (7.3)] proves that $\beta_n$ is an isomorphism for $1 \leq n < \infty$ (cf. also [3, (19.21)]). When $S'$ and $k'$ are fields of characteristic $p > 0$ André [4] proves that $\beta_n$ is an isomorphism for $1 \leq n \leq 2p$, but not necessarily for $n = 2p + 1$.

The spectral sequence in the next theorem takes an input analogous to that of Quillen's sequence described above, and is constructed along similar lines. However, it converges to invariants of $\varphi$ defined in terms of DG—rather than classical—homological algebra. We denote $\Sigma$ the suspension functor [20, (5.3)] on the category of simplicial $S$-modules.

(4.2) THEOREM. *If $\varphi \colon R \longrightarrow (S, \mathfrak{n}, \ell)$ is a surjective local homomorphism, then there is a homological first quadrant spectral sequence of $\Gamma$-algebras over $\ell$, such that*

$$^2\mathrm{E}_{p,q} = \pi_{p+q}(\mathrm{Sym}_q^\ell(\Sigma\mathcal{L}(S|R) \otimes_S \ell)) \implies \ell\langle X\rangle_{p+q}$$

*where $X = X_{\geqslant 2}$ is a set of $\Gamma$-variables with $\mathrm{card}\,(X_n) = \varepsilon_n(\varphi)$ for $n \geq 2$. The edge map*

$$\ell\langle X\rangle_n \twoheadrightarrow {}^\infty\mathrm{E}_{n-1,1} \hookrightarrow {}^2\mathrm{E}_{n-1,1} = \mathrm{D}_{n-1}(S\,|\,R,\,\ell)$$

*induces $\ell$-linear homomorphisms $\beta_n \colon \ell X_n \longrightarrow \mathrm{D}_{n-1}(S\,|\,R,\,\ell)$. If $p = \mathrm{char}\,\ell$, then $\beta_n$ is bijective for $2 \leq n < \infty$ when $p = 0$, and for $2 \leq n \leq 2p$ when $p > 0$.*

In view of Lemma (1.7), the last assertion of the theorem yields:

(4.3) COROLLARY. *For each local homomorphism $\varphi \colon R \longrightarrow (S, \mathfrak{n}, \ell)$ there are equalities*

$$\mathrm{rank}_\ell\,\mathrm{D}_n(S\,|\,R,\,\ell) = \varepsilon_{n+1}(\varphi) \ \textit{ for } \ \begin{cases} 2 \leq n < \infty & \textit{when }\,\mathrm{char}\,\ell = 0\,; \\ 2 \leq n \leq 2p - 1 & \textit{when }\,\mathrm{char}\,\ell = p > 0\,. \end{cases}$$

*Proof of the theorem.* Choose first a cofibrant simplicial algebra resolution $\mathcal{S}$ of the $R$-algebra $S$, such that $\mathcal{S}_0 = R$ (cf. [3, (9.27)]). As $\mathcal{F} = \mathcal{S}\bar{\otimes}_R\ell$ has $\mathcal{F}_0 = \ell$, choose next a cofibrant simplicial algebra resolution $\mathcal{G}$ of the $\mathcal{F}$-algebra $\ell = \mathcal{F}/\mathcal{F}_{\geqslant 1}$ by adjoining to $\mathcal{F}$ variables of degrees $\geq 2$ (cf. [3, (9.19)]). The simplicial ideal $\mathcal{J} = \mathrm{Ker}\,(\mathcal{G}\bar{\otimes}_\mathcal{F}\ell \to \ell)$ is trivial in degrees $\leq 1$, so $\pi_i(\mathcal{J}^q) = 0$ for $i < q$ by a theorem of Quillen ([37, (6.12)] or [3, (13.3)]). Thus, the $\mathcal{J}$-adic filtration of $\mathcal{G}\bar{\otimes}_\mathcal{F}\ell$ yields a convergent spectral sequence

$$^2\mathrm{E}_{p,q} = \pi_{p+q}(\mathcal{J}^q/\mathcal{J}^{q+1}) \implies \pi_{p+q}(\mathcal{G}\bar{\otimes}_\mathcal{F}\ell)\,.$$

It remains to express the vector spaces above in terms of invariants of $\varphi$.

Each $\mathcal{G}_n$ is a polynomial ring over $\mathcal{F}_n$, so there are isomorphisms of simplicial vector spaces $\mathcal{J}^q/\mathcal{J}^{q+1} \cong \mathrm{Sym}_q^\ell(\mathcal{J}/\mathcal{J}^2)$ and $\mathcal{J}/\mathcal{J}^2 = \mathcal{L}(\ell|\mathcal{F})$. The homomorphisms $\ell \to \mathcal{F} \to \ell$ give rise to a Jacobi-Zariski distinguished triangle of simplicial $\ell$-vector spaces

$$\mathcal{L}(\mathcal{F}|\ell)\bar{\otimes}_\mathcal{F}\ell \to \mathcal{L}(\ell|\ell) \to \mathcal{L}(\ell|\mathcal{F}) \to \Sigma\mathcal{L}(\mathcal{F}|\ell)\bar{\otimes}_\mathcal{F}\ell\,.$$

As $\mathcal{L}(\ell|\ell) = 0$, it produces a weak equivalence $\mathcal{L}(\ell|\mathcal{F}) \to \Sigma\mathcal{L}(\mathcal{F}|\ell)\bar{\otimes}_\mathcal{F}\ell$. On the other hand, recall that the simplicial module $\mathcal{L}(-|-)$ is defined by degreewise application of the functor $\Omega(-|-)$ of Kähler differentials. Standard change of rings properties of modules of differentials yield the isomorphism below

$$\mathcal{L}(\mathcal{F}|\ell)\bar{\otimes}_\mathcal{F}\ell = \Omega(\mathcal{F}|\ell)\bar{\otimes}_\mathcal{F}\ell \cong \Omega(\mathcal{S}|R)\bar{\otimes}_\mathcal{S}\ell = \mathcal{L}(\mathcal{S}|R)\bar{\otimes}_\mathcal{S}\ell.$$

The identification of the second page is complete.

Let $R[Y]$ be a minimal model of $S$ over $R$, as in Section 3. Since $R[Y]^\natural$ is a free graded $R$-algebra, the identity map of $S$ lifts to a morphism $R[Y] \to \mathrm{N}\mathcal{S}$ of DG algebras over $R$, that is obviously a quasi-isomorphism. For each $n$ the $R$-module $(\mathrm{N}\mathcal{S})_n$ is a direct summand of the free $R$-module $\mathcal{S}_n$, so is itself free as $R$ is local. It follows from (2.1.2) that $R[Y] \otimes_R \ell \to (\mathrm{N}\mathcal{S}) \otimes_R \ell$ is a quasi-isomorphism.

Let $\eta\colon \ell[Y] \to F$ be the composition of morphisms of DG algebras

$$\ell[Y] = R[Y] \otimes_R \ell \to (\mathrm{N}\mathcal{S}) \otimes_R \ell = \mathrm{N}(\mathcal{S} \otimes_R \ell) = \mathrm{N}\mathcal{F} = F$$

and let $\ell[Y]\langle X\rangle$ be the acyclic closure of $\ell$ over $\ell[Y]$ given by Lemma (3.7). As $G = \mathrm{N}\mathcal{G}$ is a DG $\Gamma$-algebra and $\mathrm{H}_0(G) = \ell$, a standard argument by induction on the degree of $x \in X$ shows that $\eta$ extends to a morphism of DG $\Gamma$-algebras $\zeta\colon \ell[Y]\langle X\rangle \to G$ that commutes with the divided powers of the $\Gamma$-variables in $X$; it is necessarily a quasi-isomorphism.

Define a morphism of DG algebras $\chi\colon F\langle X\rangle \to G$ by $\chi(f \otimes 1) = \eta(f)$ and $\chi(1 \otimes x^{(i)}) = \zeta(x)^{(i)}$. Note that $\chi \circ (\eta \otimes_{\ell[Y]} \ell[Y]\langle X\rangle) = \zeta$. Since $\ell[Y]\langle X\rangle^\natural$ is free over $\ell[Y]^\natural$, the map

$$\ell[Y]\langle X\rangle = \ell[Y] \otimes_{\ell[Y]} \ell[Y]\langle X\rangle \xrightarrow{\eta \otimes_{\ell[Y]} \ell[Y]\langle X\rangle} F \otimes_{\ell[Y]} \ell[Y]\langle X\rangle$$
$$= F \otimes_{\ell[Y]} \ell[Y]\langle X\rangle = F\langle X\rangle$$

is a quasi-isomorphism by (2.1.1). Thus, $\chi$ is a quasi-isomorphism.

The $F^\natural$-module $F\langle X\rangle^\natural$ is free and $G$ is the normalization of the cofibrant $\mathcal{F}$-module $\mathcal{G}$, so $\chi \otimes_F \ell$ is a quasi-isomorphism by (2.2). On the other hand, the inclusion $\partial(\ell[Y]\langle X\rangle) \subseteq (Y)\ell[Y]\langle X\rangle$ provided by (3.7) yields $\partial(F\langle X\rangle) \subseteq (F_{\geqslant 1})F\langle X\rangle$. Thus,

$$\ell\langle X\rangle = \mathrm{H}(F\langle X\rangle \otimes_F \ell) \xrightarrow{\mathrm{H}(\chi \otimes_F \ell)} \mathrm{H}(G \otimes_F \ell)$$

is an isomorphism of graded $\ell$-algebras that respects the divided powers of the $\Gamma$-variables $X_{\mathrm{even}}$. They determine the divided powers of all elements of even degree of $\ell\langle X\rangle$, so this is an isomorphism of $\Gamma$-algebras over $\ell$.

The ideals $\mathrm{N}(\mathcal{J}^q)$ are closed with respect to the divided powers of $\mathrm{N}(\mathcal{G}\bar{\otimes}_\mathcal{F}\ell)$, so our spectral sequence is one of $\Gamma$-algebras. Products of elements of positive degree and nontrivial divided powers of elements of positive even degree are contained in $\mathrm{N}(\mathcal{J}^2)$. Thus, the edge map annihilates them, and so induces an $\ell$-linear homomorphisms

$$\beta_n \colon \ell X_n \longrightarrow \pi_n(\mathcal{J}/\mathcal{J}^2) = \mathrm{D}_{n-1}(S \,|\, R, \, \ell)\,.$$

When char $\ell = 0$ the spectral sequence degenerates and each $\beta_n$ is an isomorphism by the argument of Quillen for the proof of [37, (7.3)]. When char $\ell = p > 0$ the edge homomorphism $\beta_n$ is bijective for $2 \leq n \leq 2p$ by the argument of André [4]. □

Next we prove a local version of Quillen's conjecture by replaying an argument from [9].

(4.4) THEOREM. *Let $\varphi \colon R \longrightarrow (S, \mathfrak{n}, \ell)$ be a local homomorphism, such that* $\mathrm{fd}_R S < \infty$.
*If* $\mathrm{D}_n(S \,|\, R, \, \ell) = 0$ *for* $n \gg 0$ *then $\varphi$ is c.i. at $\mathfrak{n}$.*

*Proof.* By Lemma (1.7) and Remark (3.5), we may assume that $\varphi$ is onto; the vector space $\mathrm{D}_n(S \,|\, R, \, \ell)$ then has finite rank for each $n$ (cf. [3, (5.12)]), and vanishes for $n \leq 0$.

By hypothesis, $\mathrm{D}_n(S \,|\, R, \, \ell) = \pi_n(\mathcal{L}(\ell|\mathcal{F}))$ is not zero for only finitely many $n$. By the Dold-Kan equivalence, the simplicial vector space $\mathcal{L}(\ell|\mathcal{F})$ is a direct sum $\bigoplus_{j=0}^m W_j$ of simplicial vector spaces, such that $\pi(W_0) = 0$ and there exist positive integers $n_1, \ldots, n_m$ for which $\pi_{n_j}(W_j) \cong \ell$ and $\pi_n(W_j) = 0$ if $n \neq n_j$. Thus,

$$\mathrm{H}(\mathrm{Sym}^\ell W_j) \cong \mathrm{H}(\mathrm{Sym}^\ell W_0) \otimes \bigotimes_{j=1}^m \mathrm{H}(\mathrm{Sym}^\ell W_j)\,.$$

By Dold and Thom [21], $\mathrm{H}(\mathrm{Sym}^\ell W_0) \cong \ell$ and $\mathrm{H}(\mathrm{Sym}^\ell W_j) = \mathrm{H}(\mathbb{Z}, n_j, \ell)$, where $\mathrm{H}(\mathbb{Z}, n, \ell)$ denotes the homology with coefficients in $\ell$ of the Eilenberg-MacLane space $\mathrm{K}(\mathbb{Z}, n)$ with unique nontrivial homotopy group $\pi_n(\mathrm{K}(\mathbb{Z}, n)) \cong \mathbb{Z}$.

Set $f_n(t) = \sum_{i=0}^\infty \mathrm{rank}_\ell \mathrm{H}_i(\mathbb{Z}, n, \ell) t^i$. If char $\ell = 0$, then $\mathrm{H}(\mathbb{Z}, n, \ell)$ is a free skew-commutative $\ell$-algebra on a single generator of degree $n$; hence $f_n(t) = 1 + t^n$ or $f_n(t) = 1/(1-t^n)$. When char $\ell = 2$ Serre computes $\mathrm{H}(\mathbb{Z}, n, \ell)$ and proves [39] that $f_n(t)$ converges in the open unit disk. When char $\ell = p > 2$ that conclusion is obtained by Umeda [42, (2.i)], based on the computation of $\mathrm{H}(\mathbb{Z}, n, \ell)$ by Cartan [18]. (Alternative derivations of $\mathrm{H}(\mathrm{Sym}^\ell W_j)$ can be found in [35], and of convergence in [9].)

The spectral sequence of Theorem (4.2) yields coefficientwise inequalities

$$\sum_{n=2}^\infty \varepsilon_n(\varphi) t^n = \sum_{n=2}^\infty \mathrm{card}\,(X_n) t^n \preccurlyeq \sum_{n=0}^\infty \mathrm{rank}_\ell \, \ell\langle X\rangle_n t^n \preccurlyeq \prod_{j=1}^m f_{n_j}(t)$$

so $\sum_{n=0}^{\infty} \varepsilon_n(\varphi) t^n$ converges in the unit disk. Theorem (3.10) then implies that $\varphi$ is c.i. at $\mathfrak{n}$. □

(4.5) *Remark.* If $R$ is a regular local ring, then in each Cohen factorization of $\dot{\varphi}$ the ring $R'$ is regular; it follows from the definitions that $\varphi$ is c.i. at $\mathfrak{n}$ if and only the local ring $S$ is a complete intersection.

In the special case $S = k = \ell$, the implication (iii) $\implies$ (i) of the next corollary is a well-known consequence of basic change of rings properties of André-Quillen homology (cf. (1.6)); the converse, conjectured in [37, (11.7)], was proved by Gulliksen [28] when char $k = 0$ and by Avramov [9] when char $k > 0$.

(4.6) COROLLARY. *The following conditions are equivalent*:

(i) $R$ *is a complete intersection and* $\mathrm{D}_n(S \,|\, R, \ell) = 0$ *for* $n \gg 0$.
(ii) $S$ *is a complete intersection and* $\mathrm{D}_n(S \,|\, R, \ell) = 0$ *for* $n \gg 0$.
(iii) $R$ *and* $S$ *are complete intersections*.

*When they hold*, $\mathrm{D}_n(S \,|\, R, \ell) = 0$ *for* $n \geq 3$.

*Proof.* Set char $\ell = p$, and let $\eta \colon \mathbb{Z}_{(p)} \to S$ be the local homomorphism induced by the canonical map $\mathbb{Z} \to S$. All the assertions of the corollary follow from the theorem and the preceding remark, applied to the Jacobi-Zariski exact sequence

$$\cdots \to \mathrm{D}_{n+1}(S \,|\, \mathbb{Z}_{(p)}, \ell) \to \mathrm{D}_n(R \,|\, \mathbb{Z}_{(p)}, k) \otimes_k \ell \to \mathrm{D}_n(S \,|\, R, \ell)$$
$$\to \mathrm{D}_n(S \,|\, \mathbb{Z}_{(p)}, \ell) \to \cdots$$

of the homomorphisms $\mathbb{Z}_{(p)} \to R \to S$. □

We can now prove the remaining results announced in Section 1.

*Proof of Theorems* (1.3)–(1.5). Let $\varphi \colon R \to S$ be a homomorphism of noetherian rings. If $\varphi$ is l.c.i, then it is locally of finite flat dimension by (3.5), and has $\mathrm{D}_n(S \,|\, R, -) = 0$ for $n \gg 0$ by the localization property (1.9), and by Theorem (1.2): this is the *only if* part of Theorem (1.3). The converse—Quillen's conjecture—results from Theorem (4.4) by localization.

Assume that $(m-1)!$ is invertible in $S$, and that $\mathrm{D}_n(S \,|\, R, -) = 0$ for some $n \geq 2$, such that $2 \leq n \leq 2m-1$. If $\mathfrak{q}$ is a prime ideal of $S$, then $n \leq 2\,\mathrm{char}\,k(\mathfrak{q}) - 1$. Remarks (1.9) and (3.5), and Corollary (4.3) translate Theorem (3.4) on the rigidity of deviations into Theorem (1.4) on the rigidity of André-Quillen homology.

Corollary (4.6) and localization establish all the assertions of Theorem (1.5), except for the vanishing of $\mathrm{D}_n(S\,|\,R, -) = 0$ for $n \geq 3$ when $R$ and $S$ are locally complete intersections; in that case, they give $\mathrm{D}_n(S\,|\,R, k(\mathfrak{q})) = 0$ for $n \geq 3$ and each $\mathfrak{q} \in \mathrm{Spec}\,S$, and we conclude that the homology functors vanish by [3, (S.29)]. □

We finish this section with an application of the result of [7], on maps of vector spaces of indecomposables in Tor's, to the proof of a vanishing theorem for connecting homomorphisms in some Jacobi-Zariski exact sequences. For surjective ring homomorphisms and $i = 1$ the theorem below may be read off the proof of Theorem 3 in Rodicio [38].

(4.7) THEOREM. *Let* $\psi\colon Q \to (R, \mathfrak{m}, k)$ *and* $\varphi\colon (R, \mathfrak{m}, k) \to (S, \mathfrak{n}, \ell)$ *be local homomorphisms. If* $\mathrm{fd}_R S < \infty$, *then for each local homomorphism* $\psi\colon Q \to R$ *the connecting map* $\eth_{2i}$ *in the Jacobi-Zariski exact sequence*

$$\ldots \to \mathrm{D}_{2i}(R\,|\,Q, \ell) \to \mathrm{D}_{2i}(S\,|\,Q, \ell) \to \mathrm{D}_{2i}(S\,|\,R, \ell)$$
$$\xrightarrow{\eth_{2i}} \mathrm{D}_{2i-1}(R\,|\,Q, \ell) \to \ldots$$

*is trivial for* $1 \leq i \leq \infty$ *if* $\mathrm{char}\,k = 0$, *and for* $1 \leq i \leq p$ *if* $\mathrm{char}\,k = p > 0$.

*Proof.* Consider the commutative diagram of local homomorphisms

$$\begin{array}{ccccc}
Q & \xrightarrow{\psi} & R & \xrightarrow{\varphi} & S \\
{\scriptstyle \psi}\downarrow & & {\scriptstyle \varphi}\downarrow & & \downarrow{\scriptstyle \varepsilon} \\
R & \xrightarrow{\varphi} & S & \xrightarrow{\varepsilon} & \ell \\
\downarrow & & \downarrow{\scriptstyle \varepsilon} & & \\
k & \xrightarrow{\xi} & \ell.
\end{array}$$

The canonical map $\mathrm{Tor}^R(\xi, \ell)\colon \mathrm{Tor}^R(k, \ell) \to \mathrm{Tor}^S(\ell, \ell)$ induced by $\varphi$ is a homomorphism of $\Gamma$-algebras, and so induces for each $n$ a natural homomorphism $\mathrm{T}_n(\varphi, \xi, \ell)$ of the graded $\ell$-vector spaces defined in Remark (4.1). For each $i \geq 1$, we thus get a commutative diagram of linear maps of $\ell$-vector spaces

$$\begin{array}{ccccccc}
 & & & & \mathrm{D}_{2i}(S\,|\,R, \ell) & \xrightarrow{\eth_{2i}} & \mathrm{D}_{2i-1}(R\,|\,Q, \ell) \\
 & & & & \downarrow{\scriptstyle \mathrm{D}_{2i}(\varepsilon\,|\,R, \ell)} & & \parallel \\
\mathrm{T}_{2i}(R, k, \ell) & \xrightarrow{\beta_{2i}} & \mathrm{D}_{2i}(k\,|\,R, \ell) & \xrightarrow{\mathrm{D}_{2i}(\xi\,|\,R, \ell)} & \mathrm{D}_{2i}(\ell\,|\,R, \ell) & \longrightarrow & \mathrm{D}_{2i-1}(R\,|\,Q, \ell) \\
{\scriptstyle \mathrm{T}_{2i}(\varphi, \xi, \ell)}\downarrow & & {\scriptstyle \mathrm{D}_{2i}(\xi\,|\,R, \ell)}\downarrow & & \downarrow{\scriptstyle \mathrm{D}_{2i}(\ell\,|\,\varphi, \ell)} & & \\
\mathrm{T}_{2i}(S, \ell, \ell) & \xrightarrow{\beta_{2i}} & \mathrm{D}_{2i}(\ell\,|\,S, \ell) & = & \mathrm{D}_{2i}(\ell\,|\,S, \ell) & &
\end{array}$$

in which the third column is a segment of a Jacobi-Zariski exact sequence associated with the homomorphisms $R \longrightarrow S \longrightarrow \ell$. To describe the maps in the diagram we set char $k = p$.

The map $T_{2i}(\varphi, \xi, \ell)$ is injective for $i \geq 1$ by the main result of [7].

The map $\beta_{2i}$ is bijective for all $i$ if $p = 0$, for $1 \leq i \leq p$ if $p > 0$ (cf. (4.1)).

The map $D_{2i}(\xi \,|\, R, \ell)$ is bijective for $i \geq 1$ because in the Jacobi-Zariski exact sequence associated with the homomorphisms $R \longrightarrow k \longrightarrow \ell$ the spaces $D_n(k \,|\, \ell, \ell)$ vanish for $n \geq 2$ by (1.6.3).

From the commutativity of the diagram we conclude first that $D_{2i}(\varepsilon \,|\, R, \ell) = 0$, then that $\eth_{2i} = 0$, where $1 \leq i < \infty$ if $p = 0$ and $1 \leq i \leq p$ if $p > 0$. □

## 5. L.c.i. homomorphisms

Throughout this section, $\varphi \colon R \longrightarrow S$ denotes a homomorphism of noetherian rings.

The results that follow establish the stability of the class of l.c.i. homomorphisms. They parallel those on Gorenstein homomorphisms in [11] and on Cohen-Macaulay homomorphisms in [13], with a bonus: the base change and decomposition theorems are stronger.

We start by verifying that the concept of l.c.i. homomorphism introduced in Section 1 subsumes earlier notions, whenever they are defined. From (1.9) and (4.5) we get:

(5.1) *Structure homomorphisms.* When the ring $R$ is regular (in particular, when $R = \mathbb{Z}$) the homomorphism $\varphi$ is l.c.i. if and only if the ring $S$ is l.c.i.

For a flat homomorphism the complete intersection property is easy to describe.

(5.2) *Flat homomorphisms.* When $\varphi$ is flat, it is l.c.i. if and only if all its nontrivial fibers are complete intersections.

It suffices to establish the local statement:

(5.2.1) LEMMA. *A flat local homomorphism $\varphi \colon (R, \mathfrak{m}, k) \longrightarrow (S, \mathfrak{n}, \ell)$ is c.i. at $\mathfrak{n}$ if and only if the local ring $S/\mathfrak{m}S$ is a complete intersection.*

*Proof.* The homomorphism $\varphi$ is c.i. at $\mathfrak{n}$ if and only if $D_2(S \,|\, R, \ell) = 0$, cf. (1.8). The local ring $S/\mathfrak{m}S$ is c.i. if and only if $D_3(\ell \,|\, S/\mathfrak{m}S, \ell) = 0$ (cf. [3, (6.27), (10.20)]). The test modules are isomorphic (cf. the proof of (1.7)). □

The next result applies, in particular, to $R$-algebras essentially of finite type. It can be established by elementary arguments with regular factorizations, along the lines of the proof of [11, (6.3)], but the use of cotangent homology is more expedient.

(5.3) *Smoothable homomorphisms.* Let $\varphi = \varphi'\dot\varphi$, where $\dot\varphi \colon R \longrightarrow R'$ is regular, $R'$ is noetherian, and $\varphi' \colon R' \longrightarrow S$ is surjective with kernel $\mathfrak{a}$.

If $\varphi$ is l.c.i., then $\mathfrak{a}_{\mathfrak{p}'}$ is generated by a regular sequence for each prime ideal $\mathfrak{p}'$ of $R'$.

If $\mathfrak{a}_{\mathfrak{m}'}$ is generated by a regular sequence for each maximal ideal $\mathfrak{m}'$ of $R'$, then $\varphi$ is l.c.i.

*Proof.* In view of (1.1), the Jacobi-Zariski exact sequence

$$\mathrm{D}_2(R' \,|\, R, N) \longrightarrow \mathrm{D}_2(S \,|\, R, N) \longrightarrow \mathrm{D}_2(S \,|\, R', N) \longrightarrow \mathrm{D}_1(R' \,|\, R, N)$$

yields $\mathrm{D}_2(S \,|\, R, -) \cong \mathrm{D}_2(S \,|\, R', -)$. By (1.2), $\varphi$ is l.c.i. if and only if $\mathrm{D}_2(S \,|\, R', -) = 0$. This condition can be checked locally, either over $\operatorname{Spec} S$ or over $\operatorname{Max} S$; (1.9) translates it to: $\mathrm{D}_2(S_{\mathfrak{q}} \,|\, R'_{\mathfrak{p}'}, -) = 0$ with $\mathfrak{q}$ ranging over the corresponding subset of $\operatorname{Spec}. S$ and $\mathfrak{p}' = \mathfrak{q} \cap R'$. By (1.6.1), the last equality holds if and only if $\mathfrak{a}_{\mathfrak{p}'}$ is generated by a regular sequence. $\square$

Restating part of Theorem (1.5) with the aid of Theorem (1.2), we get:

(5.4) *Homomorphisms of* l.c.i. *rings.* When the rings $R$ and $S$ are l.c.i., the homomorphism $\varphi$ is l.c.i. if and only if it is locally of finite flat dimension.

As in the case of rings, the l.c.i. property interpolates between regularity and Gorensteinness.

(5.5) *Hierarchy.* A regular homomorphism is l.c.i.
An l.c.i. homomorphism is locally Gorenstein.

*Proof.* The first assertion results from a comparison of (1.1) and (1.2).

For the second assertion, recall from [14, (3.11)] that $\varphi$ is Gorenstein at $\mathfrak{q}$ if in some Cohen factorization of $`(\varphi_{\mathfrak{q}})$ the kernel of the surjective map is a Gorenstein ideal; ideals generated by regular sequences have that property. $\square$

In the next three theorems $\psi \colon Q \longrightarrow R$ is a homomorphism of noetherian rings.

(5.6) *Composition.* If $\psi$ and $\varphi$ are l.c.i. homomorphisms then so is $\varphi\psi$.

*Proof.* Consider the exact sequence $\mathrm{D}_2(R \,|\, Q, -) \longrightarrow \mathrm{D}_2(S \,|\, Q, -) \longrightarrow \mathrm{D}_2(S \,|\, R, -)$ on the category of $S$-modules, and invoke (1.2). $\square$

As explained in the introduction, the next theorem has consequences for the functoriality of orientations of morphisms discussed in [24, (17.4.6)].

(5.7) *Decomposition.* If $\varphi$ is locally of finite flat dimension and $\varphi\psi$ is l.c.i., then $\varphi$ is l.c.i. and $\psi$ is c.i. on the image of $^a\varphi$.

Again, it suffices to establish the local statement.

(5.7.1) LEMMA. *If $\psi\colon Q \longrightarrow (R, \mathfrak{m}, k)$ and $\varphi\colon (R, \mathfrak{m}, k) \longrightarrow (S, \mathfrak{n}, \ell)$ are local homomorphisms, such that $\operatorname{fd}_R S < \infty$ and $\varphi\psi$ is c.i. at $\mathfrak{n}$, then $\psi$ is c.i. at $\mathfrak{m}$ and $\varphi$ is c.i. at $\mathfrak{n}$.*

*Proof.* We play a game of musical chairs with Proposition (1.8) and the exact sequence

$$\mathrm{D}_3(S\,|\,R,\,\ell) \longrightarrow \mathrm{D}_2(R\,|\,Q,\,\ell) \longrightarrow \mathrm{D}_2(S\,|\,Q,\,\ell) \longrightarrow \mathrm{D}_2(S\,|\,R,\,\ell) \longrightarrow 0$$

provided by Theorem (4.7). The proposition translates the hypothesis into $\mathrm{D}_2(S\,|\,Q,\,\ell) = 0$. The sequence implies $\mathrm{D}_2(S\,|\,R,\,\ell) = 0$. The proposition then yields $\mathrm{D}_3(S\,|\,Q,\,\ell) = 0$. The sequence returns $\mathrm{D}_2(R\,|\,Q,\,\ell) = 0$. The proposition concludes that $\psi$ is c.i. at $\mathfrak{m}$. □

As $^a\varphi$ is surjective when $\varphi$ is faithfully flat, the last two theorems imply:

(5.8) *Flat descent.* When $\varphi$ is faithfully flat, the composition $\varphi\psi$ is l.c.i. if and only if both homomorphisms $\varphi$ and $\psi$ have this property.

In view of (5.1), when $Q = \mathbb{Z}$ the (de)composition theorems yield statements on the transfer of the l.c.i. property of rings along $\varphi$; in the flat case, by (5.2) they show that $S$ is l.c.i. if and only if $R$ and the nontrivial fibers of $\varphi$ are l.c.i., as proved in [6].

(5.9) *Ascent.* If $R$ and $\varphi$ are l.c.i., then so is $S$.

(5.10) *Descent.* If $S$ is l.c.i. and $\varphi$ is locally of finite flat dimension, then $R$ is l.c.i. on the image of $^a\varphi$ and $\varphi$ is l.c.i.

The functorial properties of cotangent homology give a very satisfactory result on

(5.11) *Flat base change.* Let $R' \xleftarrow{\rho} R \xrightarrow{\varphi} S$ be homomorphisms of noetherian rings, such that $S' = S \otimes_R R'$ is noetherian, and let $\varphi' = \varphi \otimes_R R' \colon R' \longrightarrow S'$ be the induced map.

(1) If $\operatorname{Tor}_n^R(S, R') = 0$ for $n > 0$ and $\varphi$ is l.c.i., then $\varphi'$ is l.c.i.
(2) If $\rho$ is faithfully flat and $\varphi'$ is l.c.i., then $\varphi$ is l.c.i.

*Proof.* Under the hypothesis of (1) the functors $D_2(S'|T,-)$ and $D_2(S|R,-)$ on the category of $S'$-modules are isomorphic (cf. [3, (4.54)]). Under the hypothesis of (2) the functors $D_2(S'|T, (-\otimes_R R'))$ and $D_2(S|R,-)\otimes_R R'$ on the category of $S$-modules are isomorphic (cf. [3, (4.58)]). Now both assertions result from (1.2). □

The localization problem for the complete intersection property asks when a homomorphism $\varphi$ that is c.i. at the maximal ideals of $S$ is actually l.c.i. By (5.4) smoothable homomorphisms have this property, but localization fails in general. For instance, the $\mathfrak{m}$-adic completion map $R \longrightarrow \widehat{R}$ of a local ring $(R,\mathfrak{m})$ is always complete intersection at $\widehat{\mathfrak{m}}$, but it is l.c.i. precisely when the *formal fibers* of $R$ are l.c.i. rings. Thus, the singularities of the formal fibers of $R$ may contain obstructions to localization of c.i. homomorphisms.

When $\varphi$ is flat, Marot [33] and Tabâa [40] prove that these are the only obstructions. The next theorem extends the result to all l.c.i. homomorphisms. Applied to the structure map $\mathbb{Z} \longrightarrow S$ it shows that if $S$ is an l.c.i. ring, then so are its localizations; cf. [6].

(5.12) *Localization.* Assume that for each $\mathfrak{n} \in \mathrm{Max}\, S$ the formal fibers of the local ring $R_{\mathfrak{n} \cap R}$ are l.c.i. rings. If $\varphi$ is c.i. at each maximal ideal of $S$, then $\varphi$ is l.c.i.

In the proof of the more refined local statement below we use a special case of the result for flat homomorphisms: $R$ and $S$ are complete and $S/\mathfrak{m}S$ is regular; in [12, §3, Step 1], [13, (5.5)] this case is handled directly, with the help of Cohen factorizations.

(5.12.1) LEMMA. *Let $\varphi\colon R \longrightarrow (S,\mathfrak{n},\ell)$ be a local homomorphism that is c.i. at $\mathfrak{n}$. If the formal fibers of $R$ are l.c.i. rings, then $\varphi$ is l.c.i. and the formal fibers of $S$ are l.c.i. rings.*

*Proof.* In view of (1.8), (5.2), (1.2), we have $D_2(S|R,\ell) = 0$ and $D_2(\widehat{R}|R,-) = 0$. We want to prove that $\varphi$ is c.i. at each $\mathfrak{q} \in \mathrm{Spec}\, S$, and that the homomorphism $\sigma_{\mathfrak{q}^*}$ induced by the completion $\sigma\colon S \longrightarrow \widehat{S}$ is c.i. at each $\mathfrak{q}^* \in \mathrm{Spec}\, \widehat{S}$ lying over $\mathfrak{q}$.

Set $\mathfrak{p} = \mathfrak{q} \cap R$, let $\rho\colon R \longrightarrow \widehat{R}$ be the completion map, let $\widehat{R} \longrightarrow R' \longrightarrow \widehat{S}$ be a Cohen factorization of the homomorphism $\widehat{\varphi}\colon \widehat{R} \longrightarrow \widehat{S}$ induced by $\varphi$, and set $\mathfrak{p}^* = \mathfrak{q}^* \cap \widehat{R} \in \mathrm{Spec}\, \widehat{R}$. As $D_2(\widehat{S}|R',\ell) \cong D_2(S|R,\ell) = 0$, we see from (1.6.1) that $D_2(\widehat{S}|R',-)$ vanishes, and then from (1.9) that $D_2(\widehat{S}_{\mathfrak{q}^*}|R'_{\mathfrak{p}'}, k(\mathfrak{q}^*)) = 0$. In the Jacobi-Zariski exact sequence

$$D_2(R'_{\mathfrak{p}'}|\widehat{R}_{\mathfrak{p}^*}, k(\mathfrak{q}^*)) \longrightarrow D_2(\widehat{S}_{\mathfrak{q}^*}|\widehat{R}_{\mathfrak{p}^*}, k(\mathfrak{q}^*)) \longrightarrow D_2(\widehat{S}_{\mathfrak{q}^*}|R'_{\mathfrak{p}'}, k(\mathfrak{q}^*))$$

the first group vanishes by the special case noted before the statement of the lemma, hence $D_2(\widehat{S}_{\mathfrak{q}^*} \,|\, \widehat{R}_{\mathfrak{p}^*}, k(\mathfrak{q}^*)) = 0$. The Jacobi-Zariski exact sequence

$$D_2(\widehat{R}_{\mathfrak{p}^*} \,|\, R_{\mathfrak{p}}, k(\mathfrak{q}^*)) \longrightarrow D_2(\widehat{S}_{\mathfrak{q}^*} \,|\, R_{\mathfrak{p}}, k(\mathfrak{q}^*)) \longrightarrow D_2(\widehat{S}_{\mathfrak{q}^*} \,|\, \widehat{R}_{\mathfrak{p}^*}, k(\mathfrak{q}^*))$$

yields $D_2(\widehat{S}_{\mathfrak{q}^*} \,|\, R_{\mathfrak{p}}, k(\mathfrak{q}^*)) = 0$. By (1.7) the composition $R_{\mathfrak{p}} \longrightarrow S_{\mathfrak{q}} \longrightarrow \widehat{S}_{\mathfrak{q}^*}$ is c.i. at $\mathfrak{q}^*\widehat{S}_{\mathfrak{q}^*}$. Since the second map is flat, the desired assertions follow from (5.7). □

Finally, the proof of [11, (6.11)] shows that the localization and flat base change properties have the following consequence.

(5.13) *Completion.* Assume that for each $\mathfrak{n} \in \operatorname{Max} S$ the formal fibers of the local ring $R_{\mathfrak{n} \cap R}$ are l.c.i. rings. If $\mathfrak{a} \subset R$ and $\mathfrak{b} \subset S$ are ideals such that $\varphi(\mathfrak{a}) \subseteq \mathfrak{b} \neq S$, then the following hold for the induced homomorphism $\varphi^* \colon R^* \longrightarrow S^*$ of ideal-adic completions.

(1) If $\varphi$ is l.c.i., then so is $\varphi^*$.
(2) If $\mathfrak{a}$ is contained in the Jacobson radical of $R$ and $\varphi^*$ is l.c.i., then so is $\varphi$.

Purdue University, West Lafayette, Indiana
*E-mail*: avramov@math.purdue.edu